\documentclass{amsart}

\usepackage[foot]{amsaddr}

\usepackage{subfig}
\usepackage{graphicx}
\usepackage{amssymb,amsfonts,amsrefs}
\usepackage{hyperref}
\usepackage{mathcommands}
\usepackage{tikz-cd}
\usepackage{enumerate}
\usepackage{array}
\hypersetup{
pdftitle={A-Hilbert-Property-for-S-Class-Groups-and-the-Gross-Conjecture
},
pdfsubject={},
pdfauthor={Julian Feuerpfeil},
pdfkeywords={}
}
\usepackage{bm}
\usepackage{pgfplots}
\pgfplotsset{compat=1.18}
\usepackage[toc,page]{appendix}

\makeatletter
\newtheorem*{rep@theorem}{\rep@title}
\newcommand{\newreptheorem}[2]{%
\newenvironment{rep#1}[1]{%
 \def\rep@title{#2 \ref{##1}}%
 \begin{rep@theorem}}%
 {\end{rep@theorem}}}
\makeatother

\newtheorem{thm}{Theorem}[section]
\newreptheorem{thm}{Theorem}
\newreptheorem{cor}{Corollary}
\newreptheorem{conj}{Conjecture}
\newtheorem{cor}[thm]{Corollary}
\newtheorem{prop}[thm]{Proposition}
\newtheorem*{thm*}{Theorem}
\newtheorem{lem}[thm]{Lemma}
\newtheorem{conj}[thm]{Conjecture}
\newtheorem{quest}[thm]{Question}

\newtheorem{heu}[thm]{Heuristic}

\newtheorem{mainthm}{Theorem}

\theoremstyle{definition}
\newtheorem{defn}[thm]{Definition}
\newtheorem{defnprop}[thm]{Definition and Proposition}
\newtheorem*{defn*}{Definition}

\newtheorem{exmp}[thm]{Example}
\newtheorem{exmps}[thm]{Examples}

\theoremstyle{remark}
\newtheorem{rem}[thm]{Remark}
\newtheorem{rems}[thm]{Remarks}

\usepackage[T1]{fontenc}

\makeatletter
\def\legendre@dash#1#2{\hb@xt@#1{%
  \kern-#2\p@
  \cleaders\hbox{\kern.5\p@
    \vrule\@height.2\p@\@depth.2\p@\@width\p@
    \kern.5\p@}\hfil
  \kern-#2\p@
  }}
\def\@legendre#1#2#3#4#5{\mathopen{}\left(
  \sbox\z@{$\genfrac{}{}{0pt}{#1}{#3#4}{#3#5}$}%
  \dimen@=\wd\z@
  \kern-\p@\vcenter{\box0}\kern-\dimen@\vcenter{\legendre@dash\dimen@{#2}}\kern-\p@
  \right)\mathclose{}}
\newcommand\artinsymbol[2]{\mathchoice
  {\@legendre{0}{1}{}{#1}{#2}}
  {\@legendre{1}{.5}{\vphantom{1}}{#1}{#2}}
  {\@legendre{2}{0}{\vphantom{1}}{#1}{#2}}
  {\@legendre{3}{0}{\vphantom{1}}{#1}{#2}}
}
\def\dlegendre{\@legendre{0}{1}{}}
\def\tlegendre{\@legendre{1}{0.5}{\vphantom{1}}}
\makeatother

\usepackage[OT2,T1]{fontenc}
\DeclareSymbolFont{cyrletters}{OT2}{wncyr}{m}{n}
\DeclareMathSymbol{\Sha}{\mathalpha}{cyrletters}{"58}

\title[Hilbert 90 for S-Class Groups and the Gross--Kuz'min Conjecture]{A Hilbert 90 Property for S-Class Groups and Connections with the Gross--Kuz'min Conjecture}

\author{Julian Feuerpfeil}

\address{Dipartimento di Matematica e Applicazioni, Università di Milano-Bicocca, 20125 Milano, Italy}

\address{FEMTO-ST, Université Marie et Louis Pasteur, 25030 Besançon, France}

\email{\href{mailto:j.feuerpfeil@campus.unimib.it}{j.feuerpfeil@campus.unimib.it}}

\begin{document}
\begin{abstract} 
Let $L/K$ be a cyclic extension of number fields, and let $S$ be a finite set of places of $K$ containing the ramified and Archimedean ones. We say that $L/K$ has the $\mathbf{cl}^S$-Hilbert~90 property if, for any generator $\sigma \in \mathrm{Gal}(L/K)$, the kernel of the arithmetic norm map $\mathrm{cl}^S(L) \to \mathrm{cl}^S(K)$ coincides with $(1 - \sigma)\mathrm{cl}^S(L)$.

We establish a criterion for the $\mathbf{cl}^S$-Hilbert~90 property that depends only on arithmetic data of the base field and does not require any knowledge of the class group of $L$. We then show that, for $\ZZ_p$-extensions, the $\mathbf{cl}^{S_p}$-Hilbert~90 property at finite layers already implies the finiteness of the coinvariants of the associated Kuz'min--Tate module, providing a new finite-level criterion for an Iwasawa-theoretic property related to the Gross--Kuz'min conjecture.

This criterion can be expressed in terms of a local map closely related to Fermat quotients, making the criterion amenable to explicit computation. Motivated by extensive numerical experiments, we formulate a conjecture predicting that this criterion is satisfied for all but finitely many primes in the totally real case and present a random matrix heuristic supporting this prediction.

\noindent \textbf{Keywords. } Hilbert’s theorem 90, $S$-class groups, Gross-Kuz'min Conjecture, cohomological Mackey functors

\noindent \textbf{MSC2020} Primary: \href{https://zbmath.org/classification/?q=cc%3A11R29}{11R29}  Secondary: \href{https://zbmath.org/classification/?q=cc%3A11R37}{11R37}, \href{https://zbmath.org/classification/?q=cc%3A11R23}{11R23} \href{https://zbmath.org/classification/?q=cc%3A20J05}{20J05}

\noindent \textbf{ORCID:} \ \href{https://orcid.org/0009-0000-0148-3348}{0009-0000-0148-3348}


\end{abstract}

\maketitle


\section*{Introduction}

{
    Hilbert's Theorem 90 is arguably one of the most fundamental results in modern algebraic number theory. The version proved by Hilbert in his Zahlbericht from 1897 \cite[Satz 90]{Hilbert1897} can be phrased as follows: Let $L/K$ be a cyclic Galois extension with Galois group $\langle \sigma\rangle$. Then for $x\in L^\times$, one has $\smash{ N_{L/K}(x)=1}$ if and only if $x=\sigma(y)y^{-1}$ for some $y\in L^\times$.
    
    This theorem has been generalized in countless ways over the years and has led to numerous developments in Galois cohomology, algebraic geometry, and class field theory. The purpose of this paper is to study a variant of it for $S$-class groups of number fields and to illustrate applications in Iwasawa theory, related to the \emph{Gross--Kuz'min} (or \emph{generalized Gross}) Conjecture.

    \begin{defn*}
        Let $L/K$ be a cyclic Galois extension of number fields with Galois group $\langle \sigma \rangle$, and let $S$ be a set of places of $K$ such that $L/K$ is unramified outside $S$. We say that $L/K$ \emph{has the $\bcl^S$–Hilbert~90 property} if for every $\alpha \in \cl^S(L)$ with trivial arithmetic norm $N(\alpha) = 1$ in $\cl^S(K)$, there exists $\beta \in \cl^S(L)$ such that $\alpha = \sigma(\beta) \beta^{-1}$.
    \end{defn*}
    The $\bcl^{S_\infty}$-Hilbert~90 property has been studied in the context of capitulation phenomena and growth of ideal classes in unramified extensions by Bembom in his PhD thesis \cite{Bembom2012} under the name \emph{Furtwängler property}. Unlike in the case of Hilbert's theorem 90, this property does not hold in general. It is therefore natural to ask for criteria ensuring its validity and to understand what arithmetic information it encodes.

    A suitable structure to formalize and generalize this question is the one of \emph{cohomological Mackey functors} (also sometimes called $G$-modulations \cite{Neukirch2008}). Although they are not strictly necessary for the statement of the problems, cohomological Mackey functors provide a flexible categorical framework, simplifying several arguments. We discuss them in Section~\ref{sec:CMF}. Bley and Boltje used cohomological Mackey functors in their paper \cite{BleyBoltje2004} to study Iwasawa invariants of elliptic curves over number fields. 

    The first goal of this article is to obtain practical criteria for verifying the $\bcl^S$-Hilbert~90 property. In \cite{QuadrelliWeigel2015} Quadrelli and Weigel showed that everywhere unramified extensions in which all places in $S$ are completely split have the $\bcl^S$-Hilbert~90 property. They used it to show a more refined version of Hilbert's theorem 94 in that situation. Theorem~\ref{thm:CriterionH90RayClassGroup} generalizes this result substantially, allowing also ramified extensions. The criterion only depends on the arithmetic of the base field and the ramification behavior of the extension. In particular, it does not require any knowledge of the class group of the extension field itself. 
    
    Using it, we are able to construct large classes of extensions of arbitrarily high degree satisfying the $\bcl^S$-Hilbert~90 property (see Section~\ref{sec:ExamplesClSH90}). As a consequence each extension $\QQ(\smash{\zeta_{p^n}})/\QQ$ for an odd prime $p$ has the $\bcl^{\{p,\infty\}}$-Hilbert~90 property (see Proposition~\ref{prop:cyclotomic extension has clS H90}), raising the natural question in which infinite families this property holds uniformly.
    
    The second part of this paper studies this question for $\ZZ_p$-extensions and arrives at the following conclusion: Under mild assumptions  the $\bcl^{S\smash{_p}}$-Hilbert~90 property on the first layer of a $\ZZ_p$-extension already implies it for every level. As a consequence it controls the coinvariants in the Iwasawa module associated to the $S_p$-class groups in the $\ZZ_p$-tower. We now explain the result in greater detail.
    
    For a $\ZZ_p$-extension $K_\infty/K$ of a number field $K$ one defines the \emph{Kuz'min--Tate module} as
    \begin{align*}
        T_p(K_\infty/K)\coloneq \varprojlim_{n} \cl^{S_p}(K_n)\otimes \ZZ_p
    \end{align*}
    where $S_p$ denotes the places above $p$ and the limit is taken with respect to the norm maps. The celebrated \emph{Gross--Kuz'min Conjecture} asserts the finiteness of $T_p(K_\infty^{\rm cy}/K)_\Gamma$, where $K_\infty^{\rm cy}/K$ is the cyclotomic $\ZZ_p$-extension and $\Gamma\coloneq \Gal(K_\infty^{\rm cy}/K)$. It is essentially due to Kuz'min \cite{Kuzmin1972} and independently to Jaulent \cite{Jaulent1986}, who expressed it in terms of the finiteness of so-called \emph{logarithmic class groups} (see \cite{Jaulent1986,Jaulent1994}). Since $T_p(K_\infty/K)$ is $\Lambda$-torsion for the Iwasawa algebra $\Lambda=\ZZ_p[\![\Gamma]\!]\cong \ZZ_p[\![T]\!]$ this is equivalent to saying that the characteristic ideal of $T_p(K_\infty^{\rm cy}/K)$ is not contained in $(T)$.
    
    This conjecture has received a lot of attention throughout the years and there has been a lot of progress, though the conjecture is still open. We refer for example to \cite{Greenberg1973,Jaulent2002,Kleine2019,Maksoud2023,Maksoud2024}.

    Motivated by this conjecture we say that a $\ZZ_p$-extension $K_\infty/K$ has the \emph{finite coinvariants property (FCP)} if $T_p(K_\infty/K)_\Gamma$ is finite. It makes sense to ask, which $\ZZ_p$-extensions of a number field have the finite coinvariants property. In fact we are not aware of an example that does not satisfy this property, although it seems plausible that such an extension exists.

    For technical reasons we now restrict our attention to $\ZZ_p$-extensions, in which all the places above $p$ are either totally ramified or completely split. We denote the collection of these $\ZZ_p$-extensions by $\mathcal{F}(K,p)$. { Furthermore, for a $\ZZ_p$-extension $K_\infty/K$ we write $m_0(K_\infty/K)$ for the minimal $m$ such that $\Gal(K_\infty/K_m)$ is contained in all the open decomposition groups $D_{v}(K_\infty/K)$ of places $v$ above $p$.} Then we have 
    \begin{mainthm}
    \label{mainthm:FirstStep}
        Let $K$ be a number field and $K_\infty/K$ a $\ZZ_p$-extension with $m_0(K_\infty/K)=0$. If $K_1/K$ has the $\bcl^{S_p}$-Hilbert~90 property, then $K_\infty/K$ has the FCP. In this situation $T_p(K_\infty/K)_\Gamma\cong \cl^{S_p}(K)_{(p)}$.
    \end{mainthm}
    As pointed out to the author by Jaulent, using techniques from $p$-adic class field theory (cf., e.g., \cite{Jaulent1986,Jaulent1994}), it is possible to show that for $K_\infty\in \mathcal{F}(K,p)$ the $\cl^{S_p}$-Hilbert~90 property for all $K_n/K$ follows from the property for $K_1/K$. For the cyclotomic $\ZZ_p$-extension one can even show that the logarithmic class group $\logCl_K$ coincides with $\cl^{S_p}(K)_{(p)}$ if and only if $K_1/K$ has the $\cl^{S_p}$-Hilbert~90 property (see Proposition~\ref{prop:Hilbert 90 and logarithmic class group}).

    Next we combine Theorem~\ref{thm:CriterionH90RayClassGroup} with a theorem by Gras on the existence of cyclic extensions of degree $p$ with prescribed ramification over $S$-Hilbert class field of a number field $K$ (see Theorem~\ref{thm:GrasMunnierWeakForm} and Corollary~\ref{cor:GrasMunnierRayClassGroup}). This way one can make the Hilbert~90 property for $K_1/K$ more explicit and uniform throughout all $\ZZ_p$-extensions in $\mathcal{F}(K,p)$. For that we define for $\mathfrak{p}$ a place of $K$ above $p$ a map
    \begin{align*}
        \tau^{\mathfrak{p}}:E^{\mathfrak{p}}\to {\bigoplus}_{\mathfrak{p}\neq \mathfrak{q}\in S_p}U_{\mathfrak{q}}/U_{\mathfrak{q}}^p,\quad \alpha \mapsto {\sum}_{\mathfrak{p}\neq \mathfrak{q}\in S_p} \overline\alpha_{\mathfrak{q}}.
    \end{align*}
    {Here $E^{\fp}=\{\alpha \in K^\times:v(\alpha)=0 \text{ for }v\neq v_{\fp}\}$ denotes the set of $\fp$-units and $U_\mathfrak{
    q}$ the local units.}
    \begin{mainthm}
    \label{mainthm:TauCriterion}
        If the map $\tau^\mathfrak{p}$ is surjective for some place $\mathfrak{p}$ of $K$ above $p$, then each $K_\infty\in \mathcal{F}(K,p)$ has the FCP. If additionally to the surjectivity of $\tau^{\mathfrak{p}}$ one has $h_p(K)=1$, then each $\ZZ_p$-extension $K_\infty/K$ with $m_0(K_\infty/K)=0$ has the FCP.
    \end{mainthm}
    Again if $K_\infty^{\rm cy}\in \mathcal{F}(K,p)$, then the statement for the cyclotomic $\ZZ_p$-extension can readily be derived from techniques in $p$-adic class field theory (see Proposition~\ref{prop:tau for cyclotomic and logarithmic class group}).

    In some cases the map $\tau^{\mathfrak{p}}$ admits concrete interpretations. For quadratic fields, its surjectivity can be expressed in terms of valuations of suitable Fermat quotients see Section~\ref{sec:ImaginaryQuadraticFermat} and hence is also connected to the Wieferich prime problem. Using this, we also see under which conditions the anticyclotomic $\ZZ_p$-extension of an imaginary quadratic number fields has the FCP. In the real quadratic setting, the map $\tau^{\mathfrak{p}}$ is closely related to a type of $p$-adic regulator defined by Gras in \cite{Gras2016}, and we will see that for $p$-rational fields, the conditions of Theorem~\ref{mainthm:TauCriterion} are satisfied.
    
    Although we do not pursue this viewpoint further in the present article, there seems to be an intriguing link between the map $\tau^\mathfrak{p}$ and the arithmetic of spin symbols.

    In Section~\ref{sec:Probability} we generalize a heuristic which is usually applied to the Wieferich prime problem to the map $\tau^\mathfrak{p}$. Modeling $\tau^\mathfrak{p}$ as a random matrix over $\FF_p$ of appropriate dimensions leads to Heuristic~\ref{heu:ProbilityTauSurjective}. In Section~\ref{sec:StatisticsOnRank} we present results of numerical experiments conducted to study this behavior. For this matter, we used the OSCAR computer algebra system (see \cite{OSCAR}). The respective source code can be found on the author's GitHub \cite{Feuerpfeil2026Code}.

    The heuristic and the numerical experiments lead us to the following conjecture, which implies an earlier conjecture by Jaulent (see \cite[Conj. 11]{Jaulent2017}):
    \begin{repconj}{conj:FailiureBound}
        Let $K$ be a totally real number field. There are at most finitely many primes $p$, unramified in $K/\QQ$ for which the map $\tau^{\mathfrak{p}}$ is not surjective for all $\mathfrak{p}$ above $p$. 
    \end{repconj}
\subsection*{ Paper organization}
    Section~\ref{sec:CMF} introduces and recalls the framework of cohomological Mackey-functors. In Section~\ref{sec:GrasMunnierTheorem} we recall two variants of a theorem by Gras, which will later be necessary for the construction of examples and the proof of Theorem~\ref{mainthm:TauCriterion}. Section \ref{sec:CriterionClSH90} develops criteria for the $\bcl^S$-Hilbert~90 property, allowing us to produce several examples in Section~\ref{sec:ExamplesClSH90}. In Section \ref{sec:H90ForZpExtensions} we study how the $\bcl^{S_p}$-Hilbert~90 property behaves in $\ZZ_p$-extensions of number fields and prove Theorem~\ref{mainthm:FirstStep} and Theorem~\ref{mainthm:TauCriterion}. In Section~\ref{sec:numerics and statistics} we first derive a heuristic, asserting how often the conditions of Theorem~\ref{mainthm:TauCriterion} are satisfied, discuss implications and then present results from numerical computations.
\subsection*{Notations and conventions}

Since notation in class field theory varies throughout the literature, we collect here the conventions used throughout this article.

\begin{itemize}

\item \textbf{Class groups and units.}
We write $\cl^S(K)$ and $\cl^S_{\mathfrak m}(K)$ for the $S$-class group and the $S$-ray class group of $K$ modulo the modulus $\mathfrak m$ (i.e., a finite formal product of places of $K$), respectively. Furthermore, 
$E^S_K\coloneq \mathcal O_{K,S}^\times$,
denotes the group of $S$-units of $K$, and we write
$E_K=E_K^\emptyset=\mathcal O_K^\times$. If $K$ is clear from the context, we simply write $E^S$ instead of $E_K^S$.

\item \textbf{Artin symbols.}
If $L/K$ is Galois and $\mathfrak p$ is unramified in $L/K$ and $\mathfrak{P}$ a place of $L$ above $\mathfrak{p}$, then $(\mathfrak P,L/K)$ denotes the Frobenius automorphism at $\mathfrak P$.
Whenever the value is independent of the choice of $\mathfrak P$ (for example if $L/K$ is abelian), we simply write
$(\mathfrak p,L/K)$.

\item \textbf{Class fields.}
We denote by $H_K^S$ the $S$-Hilbert class field of $K$, and by
$R_{\mathfrak m}^S(K)$
the class field corresponding to the $S$-ray class group modulo $\mathfrak m$.

\item \textbf{Local notation.}
For a local field $K_v$, we write
$U_v=\mathcal O_{v}^\times$
for its unit group and $
U_v^{(n)}=1+\mathfrak m_v^n$
for the $n^{\rm th}$ principal unit group.
\end{itemize}
}
\section{Cohomological Mackey functors}
\label{sec:CMF}
We now define \emph{Cohomological Mackey functors}, which are special kinds of Mackey functors, satisfying an additional relation to the usual Mackey identities. Mackey functors are also known under the name $G$-\emph{modulation} (see \cite[Def. 1.5.10]{Neukirch2008}) and have been studied for a long time. They first appeared under the name \emph{cohomological Mackey functor} in \cite{Dress1973}.

Usually they are defined with respect to the full subgroup lattice of a (finite) group, but here we are only interested in the data associated with the trivial subgroup and the whole group. Thus, for us a \emph{cohomological Mackey functor} $\bX$ associated to a finite group $G$ is a diagram 
\begin{equation*}
    \begin{tikzcd}
        \bX_1\arrow[r,shift left, "N_{\bX}"]&\arrow[l,shift left,,"I_\bX"]\bX_G
    \end{tikzcd}
\end{equation*}
in the category of $G$-modules. We require that $\bX_G$ is a trivial $G$-module and for $x\in \bX_1$ and $y\in \bX_G$ we have $(I_\bX\circ N_\bX)(x)=\sum_{g\in G}gx$  and $(N_\bX\circ I_\bX)(y)=|G|\cdot y$. 
\begin{exmps}
\label{exmp:CMF}
    \begin{enumerate}
        \item \label{it:CMFCohmologyGroups} Let $H$ be a finite (or profinite group) and $N$ a normal (open) subgroup and $G\coloneq H/N$ and $M$ a (discrete) $H$-module. Group cohomology with coefficients in $M$ of $H$ and $N$ yields a cohomological Mackey functor $\mathbf{H^n}(\bl,M)$ given by the diagram
        \begin{equation*}
        \begin{tikzcd}
            \rmH^n(N,M)\arrow[r,shift left, "\mathrm{cores}"]&\arrow[l, shift left,"\mathrm{res}"]\rmH^n(H,M).
        \end{tikzcd}
        \end{equation*}
        \item \label{it:CMFab} If $H,N$ and $G$ are as in (\ref{it:CMFCohmologyGroups}), then $\mathbf{Ab}$ is defined by
        \begin{equation*}
        \begin{tikzcd}
            N^{ab}\arrow[r,shift left, "\mathrm{tr}"]&\arrow[l,shift left,,"\mathrm{inc}_*"]H^{ab}
        \end{tikzcd}
        \end{equation*}
        with $\mathrm{tr}$ being the transfer (Verlagerung) and $\mathrm{inc}_*$ the map induced by the inclusion is a cohomological Mackey functor. It has been studied intensively by Quadrelli and Weigel in \cite{QuadrelliWeigel2015}.
        \item \label{it:CMFunit} If $L/K$ is a finite Galois extension with Galois group $G=\Gal(L/K)$. Then we get a cohomological Mackey functor $\bX^\times_{L/K}$ by
        \begin{equation*}
        \begin{tikzcd}
            L^\times\arrow[r,shift left, "N_{L/K}"]&\arrow[l,shift left,"\mathrm{inc}"]K^\times
        \end{tikzcd}
        \end{equation*}
        
    \end{enumerate}
\end{exmps}
Cohomological Mackey functors for a group $G$ form a category with the natural definition of morphisms. This category turns out to be abelian and has enough projectives by a result of B.~Torrecillas and T.~Weigel \cite[\S3]{TorrecillasWeigel2013}. 

{ We now define the cohomological Mackey functor, we will mainly be interested in. 
\begin{defnprop}
    Let $L/K$ be a finite Galois extension of number fields with Galois group $G=\Gal(L/K)$. Let $S$ be a finite set of places of $K$ containing all Archimedean ones. Then the diagram
    \begin{equation*}
        \begin{tikzcd}
            \cl^S(L)\arrow[r,shift left, "N"]&\arrow[l,shift left,"I"]\cl^S(K)
        \end{tikzcd}
    \end{equation*}
    with $N$ being the arithmetic norm and $I$ the canonical inclusion defines a cohomological Mackey functor.
\end{defnprop}
\begin{proof}
    The Mackey identities already hold on the level of fractional ideals. Indeed, let $\mathfrak{p}$ be a prime ideal of $\mathcal{O}_K$. Let $e$, $f$, and $g$ denote the ramification index, the inertia degree and the number of places above $\mathfrak{p}$. If $\mathfrak{p}\mathcal{O}_L=\mathfrak{P}_1^e\cdots \mathfrak{P}_g^e$. Then 
    \begin{align*}
        (N\circ I)(\mathfrak{p})&=\prod_{i=1}^g N(\mathfrak{P}_i)^e=\prod_{i=1}^g\mathfrak{p}^{fe}=\mathfrak{p}^{[L:K]}\quad 
    \text{and} \\
        (I\circ N)(\mathfrak{P}_j)&=I(\mathfrak{p})^f=\prod_{i=1}^g\mathfrak{P}_i^{|D(\mathfrak{P}_i|\mathfrak{p})|}=N_G(\mathfrak{P}_j).
    \end{align*}
\end{proof}
If $L/K$ is unramified everywhere and all places in $S$ are completely split in $L/K$, then by class field theory $\bcl^S$ coincides with $\mathbf{Ab}$ for $\Gal(K^S_\emptyset/L)\unlhd \Gal(K_\emptyset^S/K)$, where $K^S_\emptyset$ is the maximal Galois extension of $K$ unramified everywhere and split in $S$.

Although it is not necessary for the definition of $\bcl^S$ we will later assume that $S$ additionally contains all places ramified in $L/K$ to apply some techniques from class field theory.
        
There is a canonical morphism of cohomological Mackey functors $\bX_{L/K}^\times \to  \mathbf{Div}_{L/K}^S$,
where $\mathbf{Div}^S_{L/K}$ is the quotient of all fractional ideals of $L$ resp. $K$ by the subgroup generated by $S$. Its cokernel in the category of cohomological Mackey functors coincides with $\bcl^S$.
}
\subsection{Cohomology of cohomological Mackey functors}
\label{sec:sectionCohomology}
The following \emph{section cohomology groups} of a cohomological Mackey functor were first introduced by Weigel in \cite{Weigel2007} and studied in \cite{TorrecillasWeigel2013}:

Let $G$ be a finite group and $\bX$ be a cohomological Mackey functor for $G$. Then its \emph{section cohomology groups} are defined to be
\begin{align*}
    &c_0(G,\bX)\coloneq\coker(N), &&k^0(G,\bX)\coloneq\ker(I),\\
    &c_1(G,\bX)\coloneq\ker(N)/\omega_G\cdot \bX_1, &&k^1(G,\bX)\coloneq\bX_1^G/\im(I).
\end{align*}
Here $\omega_G=\ker( \ZZ[G]\to \ZZ)$ is the augmentation ideal. We say that $G$ has the $\bX$-\emph{Hilbert~90 property} if $c_1(G,\bX)=0$. If $G$ is the Galois group of an extension $L/K$ and $\bX$ is associated to (arithmetic) data coming from this extension then we also say that $L/K$ has the $\bX$-Hilbert~90 property. Note that this notion differs slightly from the one used in \cite{QuadrelliWeigel2015}, as they require $k^0(G,\bX)=0$ to speak of the Hilbert~90 property.
    
These abelian groups are related by the following exact $6$-term sequence, where $\smash{\widehat{\rmH}}^i$ we denote the Tate cohomology groups of $G$:
\begin{equation}
\label{eq:sixtermsequence}
    \begin{tikzcd}
        0\arrow[r]&c_1(G,\bX)\arrow[r]&\smash{\widehat{\rmH}}^{-1}(G,\bX_1)\arrow[r]&k^0(G,\bX)\arrow[d]\\
        0&\arrow[l]k^1(G,\bX)& \arrow[l]\smash{\widehat{\rmH}}^0(G,\bX_1)&\arrow[l]c_0(G,\bX)
    \end{tikzcd}
\end{equation}

\begin{exmp}
\label{exmp:Hilbert90propertyCMF}
    \begin{enumerate}
        \item If $G=\langle \sigma\rangle $ is cyclic, then $G$ has the $\bX$-Hilbert~90 property if and only if $\ker( N_\bX)=(1-\sigma) \bX_1$.
        \item If $\smash{\bX_{L/K}^\times}$ is as in Example~\ref{exmp:CMF} (\ref{it:CMFunit}), and $G=\Gal(L/K)$ is cyclic, then $\smash{c_1(G,\smash{\bX_{L/K}^\times})}=0$ by Hilbert's theorem~90. Since the map $K^\times \to L^\times$ is injective we have $k^0(G,\smash{\bX^\times_{L/K}})=0$, which by \eqref{eq:sixtermsequence} shows $c_1(G,\bX^\times_{L/K})\cong \smash{\widehat{\rmH}}^{-1}(G,L^\times)$
        \item \label{it:AbHasH90} If $H,N$ and $G$ are as in Example~\ref{exmp:CMF} (\ref{it:CMFab}), then Quadrelli and Weigel showed in \cite[{Thm. A}]{QuadrelliWeigel2015} that $G$ has the $\mathbf{Ab}$-Hilbert~90 property.
        
    \end{enumerate}
\end{exmp}
{
Note that as a consequence of (\ref{it:AbHasH90}) an everywhere unramified extension $L/K$ of number fields has the $\bcl^S$-Hilbert~90 property, as long as $S$ consists of places that are completely split in the extension $L/K$. }
        
Furthermore, the group $k^0(G,\bcl^S)$ is well known as \emph{capitulation kernel} and central to class field theory. The group $k^1(G,\bcl^S)$ is sometimes called \emph{capitulation cokernel}. (Notice that some authors refer to the cokernel of $\cl^S(K)\to \cl^S(L)$ as capitulation cokernel.)

The groups $c_i(G,\bX)$ and $k^i(G,\bX)$ can be interpreted as $\Ext$-groups in an appropriate category (see \cite{TorrecillasWeigel2013}). Therefore there are associated long exact sequences for short exact sequences of cohomological Mackey functors. Moreover section cohomology groups commute with arbitrary coproducts, that is, if $\bX=\bigoplus_{j\in J}\bX_j$, then we also get
\begin{align*}
    c_i(G,\bX)={\bigoplus}_{j\in J} c_i(G,\bX_j)\qquad \text{and} \qquad k^i(G,\bX)={\bigoplus}_{j\in J} k^i(G,\bX_J)
\end{align*}
for $i=0,1$. If a cohomological Mackey functor has values in torsion abelian groups, we can decompose it and its cohomology groups into $p$-parts. More precisely, we have the following proposition:
\begin{prop}
\label{prop:torsionCMF}
    Let $\bX$ be a cohomological Mackey functor whose underlying abelian groups of $\bX_1$ and $\bX_G$ are torsion, then for each prime $p$ we get a cohomological Mackey functor $\bX_{(p)}$
    \begin{equation*}
        \begin{tikzcd}
            (\bX_1)_{(p)}\arrow[r,shift left]&\arrow[l, shift left](\bX_G)_{(p)}
        \end{tikzcd}
    \end{equation*}
    where by $A_{(p)}$ we denote the $p$-primary part of a torsion abelian group $A$. Moreover, there is a decomposition $\bX=\bigoplus_{p}\bX_{(p)}$ and 
    \begin{align*}
        c_i(G,\bX_{(p)})\cong c_i(G,\bX)_{(p)}\qquad \text{and} \qquad k^i(G,\bX_{(p)})\cong k^i(G,\bX)_{(p)}.
    \end{align*}
    for $i=0,1$. If $p$ is not a divisor of $|G|$, then $c_i(G,\bX_{(p)})=k^i(G,\bX_{(p)})=0$ for $i=0,1$.
\end{prop}

\section{Some Governing field constructions}
\label{sec:GrasMunnierTheorem}
In this part, we recall two variants of a theorem by Gras, which originated in work by Gras and Munnier \cite{GrasMunnier1998}. Since then many variants have been established, most of them by Gras himself. They are summarized in \cite[{Chap. V \S2}]{Gras2003}. The technique is based on the construction of so-called \emph{governing fields}, which capture information about classes of extensions of a given field. They allow to construct extensions with prescribed ramification.

\subsection{Existence of tamely ramified extensions of degree \texorpdfstring{$p$}{}}
Fix a number field $K$, $p$ a prime number and $T$ a set of finite places of $K$ such that $N(\mathfrak{q})\equiv 1 \pmod{p}$ for all places $\mathfrak{q}\in T$. In particular, $T$ does not contain any $p$-adic places. Now we define 
\begin{align*}
    Y\coloneq\{\alpha\in K^\times:(\alpha)=\mathfrak{a}^p\text{ for some fractional ideal }\mathfrak{a}\}.
\end{align*}
Let $\cl(K)[p]$ be the $p$-torsion part of $\cl(K)$. Then the following exact sequence shows that $Y/K^{\times p}$ is a finite group:
\begin{equation*}
    \begin{tikzcd}
        0\arrow[r]& EK^{\times p}/K^{\times p}\arrow[r]&Y/K^{\times p}\arrow[r]&\cl(K)[p]\arrow[r]&0
    \end{tikzcd}
\end{equation*}
The group $Y$ contains $K^{\times p}$ and $Y/K^{\times p}$ is sometimes also called the \lq\lq Selmer group of $K$\rq\rq\, (see \cite[Final Rem. of \S5.2.2]{Cohen2000}). The above exact sequence is similar to the one that relates the Selmer groups of elliptic curves to the Tate-Shafarevich group.

We set $K'\coloneq K(\mu_{p})$. The \emph{governing field} in this context is the Kummer extension $K'(\sqrt[p]{Y})$. The next theorem is due to Gras and Munnier \cite{GrasMunnier1998}. Many generalizations (for example, for the wild case or for higher powers of $p$) are treated in \cite[Chap. V \S2]{Gras2003}. 
\begin{thm}[Gras-Munnier]
\label{thm:GrasMunnier}
    There exists a cyclic extension $L/K$ of degree $p$, with exact set of ramification $T$ if and only if there are $a_\mathfrak{q}\in \FF_p^\times$ such that
    \begin{align*}
        \prod_{v\in T} {\big({\tilde{\mathfrak{q}}},{K'(\sqrt[p]{Y})/K'}}\big)^{a_\mathfrak{q}}=1.
    \end{align*}
    Here the $\tilde{\mathfrak{q}}$ are arbitrary places above $\mathfrak{q}$ in $K'$ for $\mathfrak{q}\in T$.
\end{thm}
\subsection{Extensions of degree \texorpdfstring{$p$}{} of the $S$-Hilbert class field}
\label{sec:GrasMunnierWeak}
There is also a \lq\lq weak form\rq\rq \, of the Gras theorem, utilizing another governing field construction, which we will use in Proposition~\ref{prop:infinitelyManyExamples} to verify the Hilbert~90 property of an extension constructed by Theorem~\ref{thm:GrasMunnier}. Moreover, we will apply this method in Section~\ref{ssec:construction of tau} to prove Theorem~\ref{mainthm:TauCriterion}. We only state the theorem in a special case, since it will be sufficient for our applications. For the statement in full generality, see \cite[{Chap. V \S2}]{Gras2003}.

Again, let $K$ be a number field. We denote by $p$ a prime and by $S$ and $T$ disjoint finite sets of places, where $S$ contains all Archimedean places and such that for all $\mathfrak{p}\in T$, which do not divide $p$, we have $N(\mathfrak{p})\equiv 1\pmod p$. For $\mathfrak{p}\in T$, we denote
\begin{align*}
    V_{\mathfrak{p}}^S\coloneq\{\varepsilon\in E^S: \varepsilon_{\mathfrak{p}}\in (U_{\mathfrak{p}})^{p}\}\quad \text{and}\quad \delta_{\mathfrak{p}}^S=\Gal\Bigl(K'\bigl(\sqrt[p]{E^S}\bigr)/K'\Bigl(\sqrt[p]{V\smash{_{\mathfrak{p}}^S}}\Bigr)\Bigr).
\end{align*}
For $\mathfrak{p}\in T$, let $\eta_{\mathfrak{p}}$ be the natural map $E^S\to U_{\mathfrak{p}}/U_{\mathfrak{p}}^p$. Then it is not hard to see via Kummer duality that $\smash{\delta_{\mathfrak{p}}^S\cong (\im \eta_{\mathfrak{p}})^\vee}$, where we write $A^\vee$ for the Pontryagin dual of $A$.
\begin{defn}
\label{def:AdmissibleGras}
We say that $\sigma\in \delta_{\mathfrak{p}}^S$ is \emph{admissible} in the following cases:
\begin{enumerate}
    \item $\mathfrak{p}\mid p$ and $\eta_{\mathfrak{p}}$ is not surjective (in this case any $\sigma$ is admissible);
    \item $\mathfrak{p}\mid p$ and $\sigma\neq 1$ if $\eta_{\mathfrak{p}}$ is surjective;
    \item $\mathfrak{p}\nmid p$ and there exists an $a_{\mathfrak{p}}\in \FF_p^\times$ such that for a $\tilde{\mathfrak{p}}$ in $K'$ above $\mathfrak{p}$, such that
        $\sigma=\smash{\big(\tilde{\mathfrak{p}},K'(\sqrt[p]{E^S})/K'\big)^{a_{\mathfrak{p}}}}.$
\end{enumerate}
\end{defn}
Now we can state the \lq\lq weak form\rq\rq \, of the Gras' theorem.
\begin{thm}[Gras theorem -- weak form]
\label{thm:GrasMunnierWeakForm}
    Let $K$, $S$, $T$ and $p$ be as above. Denote by $H^S$ the $S$-Hilbert class field of $K$. Then there exists a number field $L$ such that $L/H^S$ is cyclic of degree $p$, $L/K$ is abelian, and
    \begin{enumerate}
        \item $L/K$ is completely split at every place in $S$ and
        \item $L/H^S$ is totally ramified at all places of $H^S$ above $T$
    \end{enumerate}
    if and only if for all $\mathfrak{p}\in T$ there exists $\sigma_{\mathfrak{p}}\in \delta_{\mathfrak{p}}^S$ admissible such that $\prod_{v\in T}\sigma_{\mathfrak{p}}=1$.
\end{thm}
It is straightforward to deduce the following corollary.
\begin{cor}
\label{cor:GrasMunnierRayClassGroup}
    Let $S$, $T$, $p$ be as in Theorem~\ref{thm:GrasMunnierWeakForm}. Then the following statements are equivalent:
    \begin{enumerate}
        \item \label{it:TsubsetCondition} If $T'\subseteq T$ and $(\sigma_{\mathfrak{p}})\in \prod_{v\in T'}\delta_{\mathfrak{p}}^S$ is an admissible family, then $\prod_{v\in T'}\sigma_{\mathfrak{p}}\neq 1$;
        \item \label{it:EqualRayClassGroup} For any modulus $\mathfrak{m}$ supported in $T$, we have $\cl^S_\mathfrak{m}(K)_{(p)}\cong \cl^S(K)_{(p)}$; 
    \end{enumerate}    
\end{cor}
Corollary~\ref{cor:GrasMunnierRayClassGroup} is especially easy to apply if $T$ contains exactly one place $\mathfrak{p}$. In this situation condition \eqref{it:TsubsetCondition} is equivalent to showing that every admissible $\sigma\in \delta\smash{_{\mathfrak{p}}^S}$ is non-trivial. If $\mathfrak{p}\mid p$, then this is equivalent to $\smash{\eta_\fp}$ being surjective. If $N(\fp)\equiv 1\pmod p$, it is equivalent to saying that $\fp$ is not completely split in $K'(\sqrt[p]{E\smash{^S}})/K'$.

If $T$ contains only places above $p$, then Condition~(\ref{it:TsubsetCondition}) is equivalent to the condition that for any family $(\sigma_{\mathfrak{p}})_{\mathfrak{p}\in T}$ with at least one of the $\sigma_{\mathfrak{p}}$ is admissible $\prod_{v\in T}\sigma_{\mathfrak{p}}\neq 1$. Hence it is necessary that $\eta_\mathfrak{p}$ is surjective for every $\fp\in T$ and that
\begin{align*}
    \bigoplus_{\fp\in T}\delta_\fp^S\cong \langle \delta_\fp^S:\fp\in T\rangle.
\end{align*}

\section{General results on the \texorpdfstring{$\bcl^S$}{}-Hilbert~90 property}
\label{sec:GeneralH90}
We first relate the $\bcl^S$-Hilbert~90 property to local properties of the $S$-class group of $K$. The author thanks Georges Gras for bringing this class-field-theoretic approach to his attention. It led to a significant simplification of the previous proof of the next result. { It can be seen as an application of the fundamental exact sequence of \emph{genus theory} or the ambiguous class number formula.}  The previous proof using more geometric methods can be found in Appendix~\ref{app:GeometricProof}.
\begin{thm}
    \label{thm:Hilbert90CriterionViaNorms}
    The extension $L/K$ has the $\bcl^S_{(p)}$-Hilbert~90 property if and only if the following sequence is exact 
    \begin{align*}
        E^S \longrightarrow\bigoplus_{v\in S} (K^\times_v/N(L_v^\times))_{(p)}\overset{\sum_{v\in S}\inv_v}\longrightarrow\QQ_p/\ZZ_p
    \end{align*}
    where the first map is induced by $x\mapsto \sum_{v\in S}x_v$. In particular, if there is a place $v_0\in S$ such that the map $E^S\to \bigoplus_{v\in S\setminus \{v_0\}}(K_v^\times/N(L_v^\times))_{(p)}$ is surjective, then $L/K$ has the $\bcl^S_{(p)}$-Hilbert~90 property.
\end{thm}
\begin{proof}
    The central ingredient is a variant of the classical Chevalley--Herbrand Theorem or ambiguous class number formula. The following version can be found in \cite[{Thm. III.1.9}]{Jaulent1986}:
    \begin{align*}
        |\cl^S(L)^G|=|\cl^S(K)|\frac{\prod_{v\in S}n_v}{[L\!:\! K] [E^S:(E^S\cap N(L^\times))]}
    \end{align*}
    Since $G$ is cyclic, we have $|\cl^S(L)^G|=|\cl^S(L)_G|$. To compute $|c_1(L/K,\bcl^S)|$ we have to determine $|N(\cl^S(L))|$. If $H^S_K$ denotes the $S$-Hilbert class field of $K$, then we have 
    \begin{align*}
        |N(\cl^S(L))|=[H^S_K\!:\!L\cap H^S_K]=\frac{|\cl^S(K)|[L\!:\!H^S_K\cap L]}{[L\!:\! K]}.
    \end{align*}
    Since $H^S_K$ is the maximal abelian extension of $K$, in which all places in $S$ are split, it follows that
    \begin{align*}
        \Gal(L/H^S_K\cap L)\cong \langle D(v):v\in S\rangle
    \end{align*}
    whose order is the least common multiple of $n_v$ for all $v\in S$, which we denote by $n_S$. 
    
    Putting all those quantities together yields
    \begin{align*}
        |c_1(L/K,\bcl^S)|=\frac{\prod_{v\in S}n_v}{n_S\cdot [E^S\!:\!(E^S\cap N(L^\times))]}.
    \end{align*}
    By the Hasse norm theorem, the kernel of the map $E^S\to \bigoplus_{v\in S}K^\times_v/N_v(L_v^\times)$ is $E^S\cap N(L^\times)$. Furthermore, the order of the image of the sum of the invariant maps is $n_S$. Now, the statement about the $p$-parts follows from Proposition~\ref{prop:torsionCMF} together with a comparison of the prime factorizations on both sides of the last equation.
\end{proof}
    { Since   the kernel of 
    \begin{align*}
        {\bigoplus}_{v\in S} (K^\times_v/N_v(L_v^\times))_{(p)}\overset{\sum_{v\in S}\inv_v}\longrightarrow\QQ_p/\ZZ_p.
    \end{align*}
    coincides with the $p$-part of the relative étale cohomological Brauer group of the Galois cover $\mathcal{O}_{L,S}/\mathcal{O}_{K,S}$ (see \cite[Lemma A.1]{Aviles2004} or Lemma~\ref{lem:SesForBrauer}), we write $\Br(\mathcal{O}_{L,S}/\mathcal{O}_{K,S})_{(p)}$ instead of $\ker(\sum_{v\in S}\inv_v)$. Hence the $\smash{\bcl^S_{(p)}}$-Hilbert~90 property is equivalent to the surjectivity of the canonical map $E^S\to \Br(\mathcal{O}_{L,S}/\mathcal{O}_{K,S})_{(p)}$. 
}

The first application of this proposition is the following result, which allows to transfer the $\bcl^S_{(p)}$-Hilbert~90 property between two different extensions.
\begin{prop}
    \label{prop:H90Transfer}
    Let $L/K$ and $\widetilde{K}/K$ be two extensions of a number field $K$, the second one not necessarily Galois.

    Assume that $L/K$ is cyclic and has the $\smash{\bcl^S_{(p)}}$-Hilbert~90 property for some prime $p$. Furthermore, assume that $p$ is coprime to $[\widetilde{K}\!:\!K]$ and for every place $v$ of $S$, such that $p$ divides $[L_v\!:\!K_v]$, there is exactly one place $\widetilde{v}$ in $\widetilde{K}$ above $v$.

    Then the compositum $L\widetilde{K}/\widetilde{K}$ has the $\bcl_{(p)}^S$-Hilbert~90 property as well.
\end{prop}
\begin{proof}
    The assumption that $p$ and $d\coloneq [\widetilde{K}\!:\!K]$ are coprime guarantees, that $\widetilde{L}/\widetilde{K}$ with $\widetilde{L}\coloneq L\widetilde{K}$ is a Galois extension, and $\Gal(\widetilde{L}/\widetilde{K})_{(p)}\cong \Gal(L/K)_{(p)}$. The statement is trivial if $p$ does not divide $[L\!:\!K]$ by Proposition~\ref{prop:torsionCMF}, thus we may assume $p\mid [L\!:\!K]$. For any place $v\in S$, we also denote by $v$ a place in the respective extensions above $v$. By local class field theory, the following diagram is commutative:
    \begin{equation*}
        \begin{tikzcd}
            \widetilde{K}_v^\times/ N (\widetilde{L}_v^\times)\arrow[r,"\sim"]\arrow[d,"N_{\widetilde{K}_v/K_v}"] &\Gal(\widetilde{L}_v/\widetilde{K}_v)\arrow[d,"\mathrm{res}"]\\
            K_v^\times/N(L_v^\times)\arrow[r,"\sim"]&\Gal(L_v/K_v)
        \end{tikzcd}
    \end{equation*}
    Since $\stilde{L}_v=\stilde{K}_vK_v$ and the degree $[\stilde{K}_v\!:\!K_v]$ is a divisor of $d$ and hence coprime to $p$, $N_{\stilde{K}_v,K_v}$ induces an isomorphism on the $p$-parts. 
    
    Let $\smash{\widetilde{S}}$ be the set of places of $\widetilde{K}$ above $S$. Let $\alpha:\Br(\mathcal{O}_{\stilde{L},S}/\mathcal{O}_{\stilde{K},S})_{(p)}\to \Br(\mathcal{O}_{{L},S}/\mathcal{O}_{{K},S})_{(p)}$ be the map induced by $\oplus_{v\in S}N_{\tilde{K}_{v}/K_v}$. Since all the maps $N_{\stilde{K}_v/K_v}$ are isomorphisms and multiplication by $d$ is an automorphism of $\QQ_p/\ZZ_p$, we see that $\alpha$ is an isomorphism. Now consider the following commutative diagram:
    \begin{equation*}
        \begin{tikzcd}
            E^{\tilde{S}}_{\tilde{K}} \arrow[d,"N_{\widetilde{K}/K}"] \arrow[r,"\stilde{\beta}"]& \Br(\mathcal{O}_{\stilde{L},S}/\mathcal{O}_{\stilde{K},S})_{(p)}\arrow[d,"\alpha", "\sim" { rotate=90,anchor =south}]\\
            E^S_{K}\arrow[r,"\beta"]&\Br(\mathcal{O}_{{L},S}/\mathcal{O}_{{K},S})_{(p)}
        \end{tikzcd}
    \end{equation*}
    By Theorem~\ref{thm:Hilbert90CriterionViaNorms} we see that $\beta$ is surjective. Since the natural map
    \begin{align*}
        E^S_K/(E_K^S)^{d}\to \coker(N_{\tilde{K}/K})
    \end{align*}
    is surjective and the former is $d$-torsion one sees $\coker(N_{\tilde{K}/K})_{(p)}=0$. Hence by the kernel-cokernel sequence the natural map $0=\coker(\beta)\to \coker(\stilde{\beta})$ is surjective. This finishes the proof.
\end{proof}

\subsection{A criterion for the \texorpdfstring{$\bcl^S$}{}-Hilbert~90 property}
\label{sec:CriterionClSH90}
{ We now derive a criterion for the $\bcl^S$-Hilbert~90 property, which generalizes the one implied by the work of Quadrelli and Weigel in \cite{QuadrelliWeigel2015}. In order to also treat the ramified case, we first split up the set $S$ to separate the contributions coming from places in $S$, which are split, ramified and have inertia. In the proof of Theorem~\ref{thm:CriterionH90RayClassGroup} we will see that these cases can be treated mostly independently.} Let $K$ be a number field, $S$ a set of places of $K$ and $p$ a prime. We define the following subsets of $S$:
\begin{align*}
    S_s=\{v\in S:p\nmid n_v\},\quad S_f=\{v\in S:p\mid f_v\}, \quad  S_e=\{v\in S: p\mid e_v\}.
\end{align*}
Notice that the union of $S_s$, $S_f$, and $S_e$ is $S$ and that $S_s$ and $S_e\cup S_f$ are disjoint.
\begin{thm}
    \label{thm:CriterionH90RayClassGroup}
    Let $v_0\in S\setminus S_\infty$. Set $S_s'\coloneq S_s\cup \{v_0\}$, $S_f'\coloneq S_f\setminus \{v_0\}$, $S_e'\coloneq S_e\setminus \{v_0\}$, and $\mathfrak{f}'\coloneq\smash{\mathfrak{f}/v_0^{\mathfrak{f}(v_0)}}$, where $\mathfrak{f}$ is the conductor of $L/K$. If the following two conditions are satisfied, then $L/K$ has the $\smash{\bcl^S_{(p)}}$-Hilbert~90 property:
    \begin{enumerate}
        \item \label{it:ThmH90criterionClassGroup} the natural surjection $\cl^{S_s'}(K)_{(p)}\to \cl^{S_s'\cup S_f'}(K)_{(p)}$ is an isomorphism;
        \item \label{it:ThmH90criterionRayClassGroup} the natural surjection $\cl^{S_s'}_{\mathfrak{f}'}(K)_{(p)}\to \cl^{S_s'}(K)_{(p)}$ is an isomorphism;
    \end{enumerate}
\end{thm}
\begin{proof}
    The proof will be carried out in three steps. { The first one isolates the contribution coming from ramification and inertia, while the second and third steps treat those parts separately.}
    
    \textbf{Step 1:} By Theorem~\ref{thm:Hilbert90CriterionViaNorms} it is sufficient to show that the natural map
    \begin{align*}
        E^S\to \bigoplus_{v\in S\setminus \{v_0\}}(K_v^\times/N_v(L_v^\times))_{(p)},
    \end{align*}
    which we denote by $\psi$, is surjective. For $v\in S\setminus \{v_0\}$, we apply the snake lemma to the exact sequence $0\to \mathcal{O}_{v}^\times \to K_v^\times \to \ZZ\to 0$ and the local norm maps. By taking $p$-parts we get the following exact sequence:
    \begin{equation*}
        \begin{tikzcd}[column sep=small]
            0\arrow[r]&(\mathcal{O}_{v}^\times /N_v (\mathcal{O}_{v}^\times ))_{(p)}\arrow[r]&(K_v^\times /N(L_v^\times ))_{(p)}\arrow[r]&\arrow[r](\ZZ/f_v\ZZ)_{(p)}&0
        \end{tikzcd}
    \end{equation*}
    If $p\nmid f_v$ (i.e. $v\not \in S_f'$), then the right term is vanishing and if $p\nmid e_v$ (i.e. $v\not\in S_e'$) then the left most term is trivial. Now we sum those exact sequences for all $v\in S$ different from $v_0$ and get the following commutative diagram with exact lower row.
     \begin{equation*}
        \begin{tikzcd}[ cells={font=\everymath\expandafter{\the\everymath\displaystyle}},column sep=small]
            &E^{S_s'} \arrow[r]\arrow[d,"\psi_0"] &E^S \arrow[d,"\psi"] \arrow[dr,"\psi_1\coloneq\sum_{v\in S_f'} \ord_v", out=0, in = 140]\\
            0\arrow[r]&\bigoplus_{v\in S_e'}(\mathcal{O}_{K_v}^\times /N_v (\mathcal{O}_{L_v}^\times ))_{(p)}\arrow[r]&\bigoplus_{v\in S\setminus \{v_0\}}(K_v^\times /N_v(L_v^\times ))_{(p)}\arrow[r]&\arrow[r]\bigoplus_{v\in S_f'}(\ZZ/f_v\ZZ)_{(p)}&0
        \end{tikzcd}
    \end{equation*}
    We immediately see that $\psi$ is surjective, if both $\psi_0$ and $\psi_1$ are surjective.

    \textbf{Step 2:} We show that assumption~(\ref{it:ThmH90criterionClassGroup}) implies the surjectivity of $\psi_1$: Fix an $v\in S_f'$. Since we have $\langle S_f'\rangle \subseteq \langle S_s'\rangle$ in $\cl(K)_{(p)}$, we get a factorization in the form
    \begin{align*}
        \mathfrak{p}_v=(\widetilde{\alpha}_v)\mathfrak{b}\prod_{w\in S_s'}\mathfrak{p}_w^{c_w}.
    \end{align*}
    Here the $c_w$ are integers and $\mathfrak{b}$ is a (fractional) ideal in the $p'$-part of $\cl(K)$, which is supported away from $S$, and $\widetilde{\alpha}_v\in K^\times $. Let $m$ be the order of the class of $\mathfrak{b}$ in $\cl(K)$, which by assumption is coprime to $p$. 
    
    Set $\alpha_v\coloneq\widetilde{\alpha}_v^m\beta$ for $\beta$ a generator of $\mathfrak{b}^m$. Then $\alpha_v$ is supported in $v$ and $S_s'$, hence it is in particular an $S$-unit and $\ord_v(\alpha_v)=m$ is coprime to $p$. Thus the image of $\alpha_v$ in the $p$-part of $\ZZ/f_v\ZZ$ is a generator of the latter. Furthermore, $\alpha_v$ maps to $0$ in $\ZZ/f_w\ZZ$ for all $v\neq w\in S_f'$.

    Repeating this procedure for every place in $S_f'$, we get the surjectivity of $\psi_1$. 

    \textbf{Step 3:} It remains to show that $\psi_0$ is surjective. We first notice, that $\psi_0$ factors through $(\mathcal{O}_K/\mathfrak{f'})^\times $ by fundamental properties of the conductor.

    By taking $p$-parts in the exact sequence of \cite[Thm. V.1.7]{Milne2020} and applying condition (\ref{it:ThmH90criterionRayClassGroup}), the $p$-part of $(\mathcal{O}_K/\mathfrak{f'})^\times $ maps to $0$ in $\smash{\cl_{\mathfrak{f}'}^{S_s'}(K)}$ and consequently $E^{S_s'} $ surjects onto it. Since $\smash{\bigoplus_{v\in S_e'}(\mathcal{O}_{K_v}^\times /N (\mathcal{O}_{L_v}^\times ))_{(p)}}$ is a quotient of the $p$-part of $(\mathcal{O}_K/\mathfrak{f'})^\times $, this finishes the proof.
\end{proof}
\begin{rems}
\label{rem:RemarkToCriterionH90RayClassGroup}
    \begin{enumerate}
        \item The statement of Theorem~\ref{thm:CriterionH90RayClassGroup} does not depend on the specific extension $L/K$ but only on $K$ and the ramification behavior of $S$ in $L/K$. In particular, it does not require any knowledge about $\cl(L)$. This allows for an efficient way to test the validity of $\bcl^S_{(p)}$-Hilbert~90 property of $L/K$.
        \item The $\bcl^S$-Hilbert~90 property, according to Proposition~\ref{prop:torsionCMF}, follows from the $\bcl^S_{(p)}$-Hilbert~90 property for all primes $p$ dividing $[L\!:\!K]$. Furthermore, the place $v_0$ in Theorem~\ref{thm:CriterionH90RayClassGroup} can be chosen independently for the different prime divisors of $|G|$.
        \item Condition (\ref{it:ThmH90criterionClassGroup}) of Theorem~\ref{thm:CriterionH90RayClassGroup} guarantees the surjectivity of the map $\psi_1$ from the proof and condition (\ref{it:ThmH90criterionRayClassGroup}) ensures that $\psi_0$ is surjective. One can also estimate the size of $\coker \psi$ in terms of the sizes of $\coker(\psi_0)$ and $\coker(\psi_1)$, which are directly linked to the sizes of
        \begin{align*}
            \ker \left(\cl^{S_s'}(K)_{(p)}\to \cl^{S_s'\cup S_f'}(K)_{(p)}\right)\quad \text{and} \quad  \ker \left(\cl_{\mathfrak{f'}}^{S_s'}(K)_{(p)}\to \cl^{S_s'}(K)_{(p)}\right). 
        \end{align*}
    \end{enumerate}
\end{rems}
\subsection{Constructing examples}
\label{sec:ExamplesClSH90}
The next three examples aim to put the previous results and definitions into perspective. The first one  shows that the $\bcl^S$-Hilbert~90 property is in fact weaker than the vanishing of Tate cohomology.
\begin{exmp}
    \label{exmp:H90ButNotCohomological}
    Consider the fields $K=\QQ(\sqrt{-105})$ and $L=\QQ(i,\sqrt{-105})$. Then it can be shown that the extension $L/K$ is unramified everywhere and we have 
    \begin{align*}
        \cl(K)\cong (\ZZ/2)^3\qquad \text{and}\qquad \cl(L)\cong \ZZ/2\times \ZZ/4.
    \end{align*}
    Furthermore, if $\tau$ is a generator of $G\coloneq\Gal(L/K)\cong C_2$, then $\tau$ acts on $\cl(L)$ by the map $(a,b)\mapsto (a,-b)$. It is not difficult to see that $\smash{\widehat{\rmH}}^{-1}(G,\cl(L))\cong (\ZZ/2)^2$. Since $S=S_\infty$, all the places in $S$ are split and thus the extension has the $\bcl^{S_\infty}$-Hilbert~90 property.
\end{exmp}
In the following example we see that not every extension has the $\bcl^S$-Hilbert~90 property.
\begin{exmp}
    \label{exmp:NotH90Property}
    Consider the extension $\QQ(\sqrt{-65})/\QQ$ and let $S=\{2,5,13,\infty\}$ be the exact set of ramification. Then an explicit computation shows $\cl^S(\QQ(\sqrt{-65}))\cong \ZZ/2$. Since the capitulation kernel is trivial, we have by the exact sequence (\ref{eq:sixtermsequence})
    \begin{align*}
        c_1(\QQ(\sqrt{-65})/\QQ,\bcl^S)=\smash{\widehat{\rmH}}^{-1}(C_2,\cl^S(L))\cong \ZZ/2.
    \end{align*}
    This shows that this extension does not have the $\bcl^S$-Hilbert~90 property.
\end{exmp}
Next we see that the conditions derived in Theorem~\ref{thm:CriterionH90RayClassGroup} are only sufficient, but not necessary in general to conclude the $\bcl^S$-Hilbert~90 property.
\begin{exmp}
\label{exmp:Hilbert90ButNotRayCondition}
    We now show that if an extension $L/K$ has the $\bcl^S$-Hilbert~90 property, then it does not necessarily have to satisfy the conditions of Theorem~\ref{thm:CriterionH90RayClassGroup}. 

    Consider the extension $\QQ(\sqrt{-78})/\QQ$ and the exact set of ramification $S=\{2,3,13,\infty\}$. A straightforward computation shows $\cl^S(\QQ(\sqrt{-78}))=0$. Consequently, the extension has the $\bcl^S$-Hilbert~90 property for trivial reasons. The conductor is $\mathfrak{f}=2^3\cdot 3\cdot 13 \cdot \infty$. Then we get 
    \begin{align*}
        \cl_{2^3\cdot 3 \cdot \infty}(\QQ)&\cong \ZZ/2\oplus \ZZ/2\oplus \ZZ/2,\\
        \cl_{2^3\cdot 13 \cdot \infty}(\QQ)&\cong \ZZ/2\oplus \ZZ/2 \oplus \ZZ/12, \text{ and} \\
        \cl_{3\cdot 13 \cdot \infty}(\QQ)&\cong \ZZ/2\oplus \ZZ/12.
    \end{align*}
    The $2$-parts of these groups are not cyclic. Hence, condition~\ref{it:ThmH90criterionRayClassGroup} of Theorem~\ref{thm:CriterionH90RayClassGroup} cannot be satisfied, as $S_s'$ contains only one element for any choice of $v_0$.
\end{exmp}
\begin{prop}    
\label{prop:cyclotomic extension has clS H90}
    Let $p$ be an odd prime, $n\in \NN$ and $S$ an arbitrary finite set of places of $\QQ$ containing $p$ and $\infty$. Then $\QQ(\zeta_{p^n})/\QQ$ has the $\bcl^S$-Hilbert~90 property.
\end{prop}
\begin{proof}
    The conductor $\mathfrak{f}$ of $\QQ(\zeta_{p^n})/\QQ$ is of the form $p^m\cdot \infty$ for some $m$. Now we choose $v_0\coloneq p\in S$. Then condition \ref{it:ThmH90criterionClassGroup} of Theorem~\ref{thm:CriterionH90RayClassGroup} is trivially satisfied. Condition \ref{it:ThmH90criterionRayClassGroup} follows from $\cl_\infty(\QQ)\cong \cl(\QQ)=0$.
\end{proof}
Next, we show that there exist a plethora of extensions, satisfying the $\bcl^S$-Hilbert~90 property. Our main ingredient will be a theorem due to Gras in its various forms, which have been introduced and discussed in Section~\ref{sec:GrasMunnierTheorem}. Here we focus primarily on the tame situation, i.e., where no place above $p$ is contained in $S$. 
\begin{prop}
    \label{prop:infinitelyManyExamplesOfDegreeL}
    Let $K$ be a number field, $p$ be a prime number, and $\mathfrak{q}_1$ a prime of $K$ such that $N(\mathfrak{q}_1)\equiv 1\pmod p$ and the class of $\mathfrak{q}_1$ is trivial in $\cl(K)_{(p)}$. Then there exist infinitely many primes $\mathfrak{q}_2$ of $K$ and associated cyclic Galois extensions $L^{\mathfrak{q}_2}/K$ of degree $p$ such that 
    \begin{enumerate}
        \item the extension $L^{\mathfrak{q}_2}/K$ is totally ramified at ${\mathfrak{q}_1}$ and ${\mathfrak{q}_2}$ and unramified at every other place;
        \item the extension $L^{\mathfrak{q}_2}/K$ has the $\bcl^{S}$-Hilbert~90 property with respect to $S=\{{\mathfrak{q}_1},\mathfrak{q}_2\}\cup S_\infty$.
    \end{enumerate}
\end{prop}
\begin{proof}
    We use Theorem~\ref{thm:GrasMunnier} to show the existence of the extension and Theorem~\ref{thm:GrasMunnierWeakForm} in conjunction with Theorem~\ref{thm:CriterionH90RayClassGroup} to deduce the $\bcl^S$-Hilbert~90 property. Let $K'$ and $Y$ be as in Theorem~\ref{thm:GrasMunnier}. For a place $\mathfrak{r}$ of $K$, we write $\tilde{\mathfrak{r}}$ for any place of $K'$ lying above $\mathfrak{r}$. Furthermore, abbreviate $E^{\mathfrak{q}_1}\coloneq E^{\{{\mathfrak{q}_1}\}}$. The relevant governing fields introduced in Section~\ref{sec:GrasMunnierTheorem} are diagrammatically depicted as follows.
        \begin{center}
            \begin{tikzpicture}
            \node (Lq) at (0,1.2) {$L^{\mathfrak{q}_2}$};
            \node (K) at (0,0) {$K$};
            \node (K1) at (1.5,.37) {$K'$};
            \node (KE) at (1.5,1.57) {$K'(\sqrt[p]{E^{\mathfrak{q}_1}})$};
            \node (KV) at (4,1) {$K'(\sqrt[p]{Y})$};
            \node (KEV) at (4,2.2) {$K'(\sqrt[p]{E^{\mathfrak{q}_1} Y})$};
            \draw (K) -- (K1) -- (KV) -- (KEV) -- (KE) -- (K1);
            \draw[dashed] (K) -- (Lq) node [midway,left] {$\ZZ/p$};
        \end{tikzpicture}
        \end{center}

    In the first step, we show that the field $K'(\sqrt[p]{E^{{\mathfrak{q}_1}}})$ is not contained in $K'(\sqrt[p]{Y})$. By Kummer theory this would be equivalent to $E^{{\mathfrak{q}_1}} K'^{\times p}\subseteq Y{K'}^{\times p}$. Since the class of ${\mathfrak{q}_1}$ in $\cl(K)_{(p)}$ is trivial there exists an element $\alpha$ of $E^{\mathfrak{q}_1}$ such that $p \nmid v_{\mathfrak{q}_1}(\alpha)$, which is therefore not contained in $Y{K'}^{\times p}$.

    Now we show that there is $\sigma\in \Gal(K'(\sqrt[p]{E^{\mathfrak{q}_1} Y})/K')$, such that 
    \begin{align*}
        \sigma|_{K'(\sqrt[p]{Y})}=\rho\coloneq\big(\tilde{\mathfrak{q}}_1,K'(\sqrt[p]{Y})/K'\big)^{-1} \quad \text{and}\quad \sigma|_{K'(\sqrt[p]{E^{\mathfrak{q}_1}})}\neq 1.
    \end{align*}
    This is certainly possible if $\rho\not \in \Gal(K'(\sqrt[p]{Y})/F)$, where $F$ is the intersection of $K'(\sqrt[p]{Y})$ and $K'(\sqrt[p]{E^{\mathfrak{q}_1}})$. In the other case this follows from the isomorphism
    \begin{align*}
        \Gal(K'(\sqrt[p]{E^{\mathfrak{q}_1}})/F)\times \Gal(K'(\sqrt[p]{Y})/F)\cong \Gal(K'(\sqrt[p]{E^{\mathfrak{q}_1} Y})/F).
    \end{align*}
    Using the Chebotarev Density Theorem, we choose a prime $\mathfrak{Q}$ of $K'(\sqrt[p]{E^{\mathfrak{q}_1} Y})$ such that $\big(\mathfrak{Q},K'(\sqrt[p]{E^{\mathfrak{q}_1} Y})/K\big)=\sigma$.
        
    Set ${\mathfrak{q}_2}$ to be the prime of $K$ and $\tilde{\mathfrak{q}}_2$ the prime of $K'$ below $\mathfrak{Q}$. Then ${\mathfrak{q}_2}$ is completely split in $K'/K$, as $\sigma\in \Gal(K'(\sqrt[p]{E^{\mathfrak{q}_1} Y})/K')$. Thus $N({\mathfrak{q}_2})\equiv 1\pmod{p}$ and we can construct the desired extension by  Theorem~\ref{thm:GrasMunnier}. By construction of $\rho$
    \begin{align*}
        \big(\tilde{\mathfrak{q}}_1,K'(\sqrt[p]{Y})/K'\big)\big(\tilde{\mathfrak{q}}_2,K'(\sqrt[p]{Y})/K'\big)=1.
    \end{align*}
    Thus the conditions of Theorem~\ref{thm:GrasMunnier} are satisfied. Denote the resulting degree $p$ extension given by the theorem by $L^{\mathfrak{q}_2}$. It is totally ramified at ${\mathfrak{q}_1}$ and ${\mathfrak{q}_2}$. 
    
    It remains to show the $\bcl^{S}$-Hilbert~90 property for $L^{\mathfrak{q}_2}/K$ and $S=\{{\mathfrak{q}_1},{\mathfrak{q}_2}\}\cup S_\infty$, which is equivalent to the $\bcl^S_{(p)}$-Hilbert~90 property. Since the ramification is tame the conductor of the extension is $\mathfrak{f}\coloneq{\mathfrak{q}_1}{\mathfrak{q}_2}$. By showing $\cl^{{\mathfrak{q}_1}}_{{\mathfrak{q}_2}}(K)_{(p)}\cong \cl^{{\mathfrak{q}_1}}(K)_{(p)}$ we deduce the $\bcl^{S}$-Hilbert~90 property from Theorem~\ref{thm:CriterionH90RayClassGroup}. But this follows from Corollary~\ref{cor:GrasMunnierRayClassGroup} by choosing $S'=\{{\mathfrak{q}_1}\}\cup S_\infty$ and $T=\{{\mathfrak{q}_2}\}$ since 
    \begin{align*}
        \big(\tilde{\mathfrak{q}}_2,K'(\sqrt[p]{E^{{\mathfrak{q}_1}}})/K'\big)=\sigma|_{K'(\sqrt[p]{E^{\mathfrak{q}_1}})}\neq 1.
    \end{align*}
\end{proof}
By estimating the extension degree $[K'(\sqrt[p]{YE^{\mathfrak{q}_1}})\!:\!K]$ it is also possible to estimate the density of primes ${\mathfrak{q}_2}$, for which the conditions are satisfied by the Chebotarev density theorem. An upper bound is given by
\begin{align*}
    h_{p}(K)p^{2(r_1+r_2)+1}
\end{align*}
with $(r_1,r_2)$ being the signature of $K$. Thus the density of the primes ${\mathfrak{q}_2}$ satisfying the conditions, is at least the reciprocal of that value. 
    
Using a more general version of Gras's theorem, it is also possible to extend this result to extensions of arbitrarily high degree. 
\begin{prop}
\label{prop:infinitelyManyExamples}
    Let $K$ be a number field, $p$ be a prime number, $r\in \NN$ such that the $p$-part of $\cl(K)$ is elementary abelian if $r>1$ and $r=1$ if $p=2$. Let ${\mathfrak{q}_1}$ be a prime of $K$ such that $N({\mathfrak{q}_1})\equiv 1 \pmod{p^r}$ and the class of ${\mathfrak{q}_1}$ is trivial in $\cl(K)_{(p)}$.

    Then there exist infinitely many primes ${\mathfrak{q}_2}$ of $K$ and associated cyclic Galois extensions $L^{{\mathfrak{q}_2}}/K$ of degree $p^r$ such that 
    \begin{enumerate}
        \item the extension $L^{\mathfrak{q}_2}/K$ is totally ramified at ${\mathfrak{q}_1}$ and ${\mathfrak{q}_2}$ and unramified at every other place;
        \item the extension $L^{\mathfrak{q}_2}/K$ has the $\bcl^{S}$-Hilbert~90 property with respect to $S=\{{\mathfrak{q}_1,\mathfrak{q}_2}\}\cup S_\infty$.
    \end{enumerate}
\end{prop}

\section{Applications to \texorpdfstring{$\ZZ_p$}{} extensions of number fields}
\label{sec:H90ForZpExtensions}
{
Up to this point we have only investigated the $\bcl^S$-Hilbert~90 property in finite cyclic extensions. Proposition~\ref{prop:cyclotomic extension has clS H90} suggests that there might be a uniform behavior in certain infinite families of extensions. 

We therefore turn to $\ZZ_p$-extensions. The key observation underlying this part of the paper is that the $\bcl^S$-Hilbert~90 property yields bounds on the growth of $\cl^{\smash{S}}(K_n)$ in a tower $K_0\subseteq K_1\subseteq K_2\subseteq ...$ of cyclic extensions and in favorable situations,  the Hilbert~90 property at a finite level can already yield by the Nakayama Lemma bounds on all levels.}

Let $K$ be a number field and $p$ be a prime and $K_\infty/K$ a $\ZZ_p$-extension of $K$ with Galois group $\Gamma\cong \ZZ_p$. We denote by $K_n/K$ the unique subextension of degree $p^n$ and define $\Gamma_n\coloneq\Gal(K_\infty/K_n)$, which is again isomorphic to $\ZZ_p$. Such an extension is only ramified at the places above $p$.

The completed group algebra $\Lambda\coloneq\ZZ_p[\![ \Gamma ]\!]=\varprojlim \ZZ_p[\Gamma/\Gamma_n]$ (see \cite[{Chap. V, \S2}]{Neukirch2008}) is a local ring of dimension $2$ and the choice of a topological generator $\gamma$ of $\Gamma$ defines an isomorphism to $\ZZ_p[\![T]\!]$, via $\gamma\mapsto T+1$.

\begin{defn}
\label{def:KuzmintateModule}
    Let $K$ be a number field and $K_\infty/K$ be a $\ZZ_p$-extension of $K$. Set $S_p$ to be the set of places of $K$ lying above $p$, together with the Archimedean ones. We call the compact $\Lambda$-module
    \begin{align*}
        T_p(K_\infty/K)\coloneq\varprojlim \cl^{S_p}(K_n)_{(p)},
    \end{align*}
    where the limit is taken with respect to the arithmetic norm maps, the \emph{Kuz'min-Tate} module associated to $K_\infty/K$.
\end{defn}
The $\Lambda$-module $T_p(K_\infty/K)$ is finitely generated, compact, and torsion (see \cite[{Prop. IX.11.1.4}]{Neukirch2008}). Moreover, it can be identified with $\Gal(H_\infty/K_\infty)$, where $H_\infty$ is the compositum of the maximal unramified abelian $p$-extension $H_n/K_n$, which are completely split at all places of $K_n$ lying above $p$. 

For a $\ZZ_p$-extension $K_\infty/K$ we denote by $m_0=m_0(K_\infty/K)$ the unique integer such that $\Gamma_{m_0}$ is equal to the intersection of all the decomposition groups of $\mathfrak{p}\in S_p$, which are open in $\Gamma$ (see \cite[Lem.~11.1.5]{Neukirch2008}). 

{ We write $\mathcal{F}(K,p)$ for the set of $\ZZ_p$-extensions $K_\infty$ over $K$ for which each place above $p$ is totally ramified or completely split. In particular for every $K_\infty\in \mathcal{F}(K,p)$ one has $m_0(K_\infty/K)=0$.}

\subsection{The Gross--Kuz'min Conjecture}
\label{sec:LogarithmicClassGroup}
{Let $\QQ_\infty^{\rm cy}$ be the subfield of $\QQ(\mu_{p^\infty})$ fixed by the torsion subgroup of $\Gal(\QQ(\mu_{p^\infty})/\QQ)\cong \ZZ_p^\times$. It is called the \emph{cyclotomic} $\ZZ_p$-extension of $\QQ$. Note that it is the only $\ZZ_p$-extension of $\QQ$. For an arbitrary number field one sets $K_\infty^{\rm cy}\coloneq K\QQ_\infty^{\rm cy}$. It is the most studied $\ZZ_p$-extension in Iwasawa theory.}

The abelian pro-$p$ group ${T_p(K_\infty^{\rm cy}/K)}_\Gamma$ then agrees with the logarithmic class group $\logCl_K$ as defined by Jaulent in \cite{Jaulent1994}. The following conjecture was first formulated by Kuz'min in this generality (see \cite{Kuzmin1972}) and independently by Jaulent (see \cite{Jaulent1986}) and is widely known as \emph{Gross--Kuz'min Conjecture} or \emph{generalized Gross Conjecture}.
\begin{conj}[Gross--Kuz'min]
    \label{conj:GrossKuzmin}
    The group ${T_p(K_\infty^{\rm cy}/K)}_\Gamma\cong \logCl_K$ is finite. 
\end{conj}
There are several equivalent formulations related to the non-vanishing of a $p$-adic regulator \cite{Gross1981} and vanishing orders of $L$-functions. Most of the methods dealing with this use $p$-adic transcendence theory. We focus on the Iwasawa theoretic approach. 

This conjecture is known to hold if $K/\QQ$ is an abelian extension by a theorem by Greenberg \cite{Greenberg1973}. It is true by Chevalley's Theorem \cite[{Lem. 4.1 Chap. 13}]{Lang1983} if there is exactly one prime in $K$ above $p$. In 2019, Kleine proved it for the case where there are two primes above $p$ (cf. \cite{Kleine2019}). Recently, Maksoud extended the list of fields for which the Gross--Kuz'min Conjecture holds substantially (cf. \cite[{Thm. 4.3.2}]{Maksoud2023}) and proved its validity for abelian extensions of imaginary quadratic fields (cf. \cite{Maksoud2024}). Important criteria, which depend on the ramification and splitting behavior of $p$ in $K/\QQ$ have been given by several authors (see for example \cite{Jaulent2002,Kuzmin2018}).

{
Since the property claimed by the Gross--Kuz'min Conjecture make sense for every $\ZZ_p$-extension, we pose the following definition.
\begin{defn}
    Let $K$ be a number field and $p$ be a prime. A $\ZZ_p$-extension $K_\infty/K$ is said to have the \emph{finite coinvariants property (FCP)} if $T_p(K_\infty/K)_\Gamma$ is finite.
\end{defn}
The Gross--Kuz'min Conjecture is equivalent to the FCP for the cyclotomic $\ZZ_p$-extension of a number field. We are not aware of an example of a $\ZZ_p$-extension of a number field that does not satisfy the FCP. 
}

\subsection{The \texorpdfstring{$\bcl^{S_p}$}{}-Hilbert~90 property implies FCP}
\label{sec:H90AndGrossKuzMin}
In the remainder of the article, we are only interested in the $\smash{\bcl^{S_p}_{(p)}}$-Hilbert~90 property, where $S_p$ denotes the places above $p$ and the infinite ones. To simplify notation, we abbreviate this cohomological Mackey functor by
\begin{align}
\label{eq:abbreviationOfClSp}
    \bcl_{[p]}\coloneq\smash{\bcl^{S_p}_{(p)}}.
\end{align}
and we write $\cl_{[p]}(K)$ for $\cl^{S_p}(K)_{(p)}$. The brackets are intended to avoid confusion with ray class groups. Note that for the extensions we consider, the $\bcl^{S_p}$-Hilbert~90 property is equivalent to the $\bcl_{[p]}$-Hilbert~90 property, as they are of $p$-power degree. However, we prefer to use $\bcl_{[p]}$, as its underlying abelian groups coincide with those of interest for us.

The following lemma is classical, and a proof for invariants and the $p$-part of class groups can be found in \cite[{Thm. 1}]{Fukuda1994}. It uses the topological Nakayama lemma and \cite[{Lem. 11.1.5}]{Neukirch2008}.
\begin{lem}
    \label{lem:StabilizationOfCoinvariants}
    If there exists $n\geq m_0(K_\infty/K)$ such that the norm map induces an isomorphism ${\cl_{[p]}(K_{n+1})}_{\Gamma}\cong {\cl_{[p]}(K_{n})}_{\Gamma}$, then ${T_p(K_\infty/K)}_\Gamma$ is finite. More precisely, we have ${T_p(K_\infty/K)}_\Gamma\cong {\cl_{[p]}(K_{n})}_{\Gamma}$.
\end{lem}
\begin{prop}
\label{prop:stabilizationOfC1}
    If there exist $m\geq n\geq m_0$ such that the norm map $N_{m+1,m}$ induces an isomorphism
    \begin{align*}
        c_1(\Gamma_n/\Gamma_m,\bcl_{[p]})\cong c_1(\Gamma_n/\Gamma_{m+1},\bcl_{[p]}),
    \end{align*}
    then $K_\infty/K$ has the FCP and $T_p(K_\infty/K)_{\Gamma_n}\cong \cl_{[p]}(K_m)_{\Gamma_n}$.

    { If $m_0(K_\infty/K)=0$ and $K_1/K$ has the $\smash{\bcl_{[p]}}$-Hilbert~90 property, then $T_p(K_\infty/K)_\Gamma\cong \cl_{[p]}(K)$ and is therefore finite. Moreover, $K_n/K$ has the $\bcl_{[p]}$-Hilbert~90 property for all $n\geq 1$.}
\end{prop}
\begin{proof}
    For the proof, we abbreviate $T_p(K_\infty/K)$ by $T$. For $m\geq n\geq m_0$, the norm maps 
    \begin{align*}
        N_{m,n}:\cl_{[p]}(K_m)\to \cl_{[p]}(K_n)
    \end{align*}
    are surjective. Hence we get the following commutative diagram with exact rows:
    \begin{equation*}
        \begin{tikzcd}
            0\arrow[r]&c_1(\Gamma_n/\Gamma_{m+1},\bcl_{[p]})\arrow[r] \arrow[d]&{\cl_{[p]}(K_{m+1})}_{\Gamma_n}\arrow[r,"N_{m+1,n}"]\arrow[d,"N_{m+1,m}"] &\cl_{[p]}(K_n)\arrow[r] \arrow[d,equal]&0\\
            0\arrow[r] &c_1(\Gamma_n/\Gamma_m,\bcl_{[p]})\arrow[r] &{\cl_{[p]}(K_{m})}_{\Gamma_n}\arrow[r,"N_{m,n}"] &\cl_{[p]}(K_n)\arrow[r] &0
        \end{tikzcd}
    \end{equation*}
    By assumption, the leftmost arrow is an isomorphism and so is the center one by the five lemma and we can apply Lemma~\ref{lem:StabilizationOfCoinvariants}, which yields $T_{\Gamma_n}\cong {\cl_{[p]}(K_{m})}_{\Gamma_n}$. Since $\cl_{[p]}(K_{m})$ is finite, so is $T_{\Gamma_n}$ and therefore $T_{\Gamma}$.

    {
    The second statement follows by setting $n=m=0$, noting that $c_1(1,\smash{\bcl_{[p]}})=0$ for trivial reasons. The isomorphism can be deduced from the surjectivity of the norm maps and the general fact that $\smash{c_1(\Gamma/\Gamma_n,\bcl_{[p]})}=0$ is equivalent to $(T_n)_\Gamma=N_{n,0}(T_n)=T_0$.}
\end{proof}
\begin{rem}
{
    It is possible to show the second statement of Proposition~\ref{prop:stabilizationOfC1} by considering the inverse limit of the sequences of Theorem~\ref{thm:Hilbert90CriterionViaNorms} for $L=K_n$. The groups in the limit sequence are finitely generated abelian pro-$p$ groups and hence by the Nakayama lemma it is possible to determine the exactness in the middle by mod-$p$ reduction. The author thanks Jaulent for pointing out this alternative proof.}
\end{rem}
{
If $K_\infty/K$ is the cyclotomic $\ZZ_p$-extension, the second statement of Proposition~\ref{prop:stabilizationOfC1} is an easy consequence of $p$-adic class field theory and can even be strengthened. This observation was communicated to the author by Jaulent in a private correspondence. The author is grateful to him for sharing the following results and its proof. 
\begin{prop}
    \label{prop:Hilbert 90 and logarithmic class group}
    Let $K_\infty/K$ be the cyclotomic $\ZZ_p$-extension of a number field $K$ and assume that $m_0=0$. Then $\logCl_K\cong \cl_{[p]}(K)$ if and only if $K_1/K$ has the $\bcl_{[p]}$-Hilbert~90 property. 
\end{prop}
\begin{proof}
    We use the standard notation from $p$-adic class field theory (see for example \cite{Jaulent1994}). Abbreviate $L=K_1$. For each $\fp\in S_p$ is a chain of epimorphisms
    \begin{align*}
        \ZZ_p\cong \mathcal{R}_{K_\fp}/\stilde{\mathcal{U}}_{K_\fp}\twoheadrightarrow \mathcal{R}_{K_\fp}/N(\mathcal{R}_{L_\fp})\twoheadrightarrow K_\fp^\times/N(L_\fp^\times)\cong \ZZ_p/p.
    \end{align*}
    where $\mathcal{R}_{K_\fp}\cong \ZZ_p\otimes K_\fp^\times$ and $\stilde{\mathcal{U}}_{K_\fp}$ the subgroup of logarithmic units. All these maps together yield a natural map $\stilde{\mathcal{D}\ell}'\to \Br(\mathcal{O}_{L,S_p}/\mathcal{O}_{K,S_p})_{(p)}$, where $\stilde{\mathcal{D}\ell}'$ is the wild part of the degree $0$ logarithmic divisors (see the proof of Scolie 3.3 from \cite{Jaulent1994}). It identifies $\Br(\mathcal{O}_{L,S_p}/\mathcal{O}_{K,S_p})_{(p)}$ with $ \stilde{\mathcal{D}\ell}'\otimes \FF_p$. From that we get a commutative square
    \begin{equation}
    \begin{tikzcd}
        \stilde{\mathcal{E}'}\arrow[r,"\stilde{{\rm div}'}"]\arrow[d,two heads]&\stilde{\mathcal{D}\ell}'\arrow[d,two heads]\\
        E^{S_p}\otimes \FF_p\arrow[r,"\beta"]&\Br(\mathcal{O}_{L,S_p}/\mathcal{O}_{K,S_p})_{(p)}
    \end{tikzcd}       
    \end{equation}
    Assume that $\logCl_K\cong \cl_{[p]}(K)$, then the short exact sequence from the proof of Scolie 3.3 in \cite{Jaulent1994} implies that $\stilde{\rm div}_K$ is surjective. Thus by the commutativity the same is true for $\beta:E^{S_p}\to \Br(\mathcal{O}_{L,S_p}/\mathcal{O}_{K,S_p})_{(p)}$ and we get the $\bcl_{[p]}$-Hilbert~90 property from Theorem~\ref{thm:Hilbert90CriterionViaNorms}.

    Now, if $L/K$ has the $\bcl_{[p]}$-Hilbert~90 property, then $\beta$ is surjective by Theorem~\ref{thm:Hilbert90CriterionViaNorms} and we get the surjectivity of $\stilde{\rm div}'$ by the Nakayama lemma using the identification $E^{S_p}\otimes \FF_p\cong \stilde{\mathcal{E}}'\otimes \FF_p$. Hence in this situation the natural map $\logCl_K\to \cl_{[p]}(K)$ is injective and using $m_0=0$ it is surjective, showing the desired statement.    
\end{proof}
}
\subsection{The map \texorpdfstring{$\tau^{\mathfrak{p}}$}{}}
\label{ssec:construction of tau}
We start with a rather technical proposition, which gives a simple criterion to verify the conditions of Corollary~\ref{cor:GrasMunnierRayClassGroup} in a special context. It is proven by Kummer duality applied to the fields constructed in Section~\ref{sec:GrasMunnierWeak} together with linear-algebraic arguments.
\begin{prop}
\label{prop:WildResidues}
    Let $K$ be a number field and $p$ a prime number. Fix a place $\mathfrak{p}$ of $K$ above $p$. If the map 
    \begin{align*}
        \tau^{\mathfrak{p}}:E^{\mathfrak{p}}\to {\bigoplus}_{\mathfrak{p}\neq \mathfrak{q}\in S_p}U_{\mathfrak{q}}/U_{\mathfrak{q}}^p,\quad \alpha \mapsto {\sum}_{\mathfrak{p}\neq \mathfrak{q}\in S_p} \overline{\alpha_{\mathfrak{q}}}
    \end{align*}
    is surjective, then for any modulus $\mathfrak{f}$ of $K$ supported in the places above $p$, which are different from $\mathfrak{p}$ we have $\cl^{\mathfrak{p}}_{\mathfrak{f}}(K)_{(p)}\cong \cl^{\mathfrak{p}}(K)_{(p)}$.
\end{prop}
\begin{proof}
    If in the following we write $\mathfrak{q}$ we mean an arbitrary prime ideal above $p$, which is different from $\mathfrak{p}$. We apply Corollary~\ref{cor:GrasMunnierRayClassGroup} to $S=\{\mathfrak{p}\}\cup S_\infty$ and $T=S_p\setminus \{\mathfrak{p}\}$. For that, we first construct the relevant governing field $K'(\sqrt[p]{E\smash{^{\mathfrak{p}}}})$, where $K'\coloneq K(\mu_p)$. Now the objects from Section~\ref{sec:GrasMunnierTheorem} read as
    \begin{align*}
        V^S_{\mathfrak{q}}=\{\varepsilon\in E^{\mathfrak{p}}:\varepsilon_{\mathfrak{q}}\in U_{\mathfrak{q}}^p\}\qquad \text{and}\qquad \delta^S_{\mathfrak{q}}= \Gal\Big(K'\big(\sqrt[p]{E^{\mathfrak{p}}}\big)/K'\big(\sqrt[p]{V\smash{_{\mathfrak{q}}}}\big)\Big).
    \end{align*}
    The components $E^{\fp}\to U_{\mathfrak{q}}/U_{\mathfrak{q}}^p$ are identified with the maps $\eta^S_\mathfrak{q}$ from Section~\ref{sec:GrasMunnierWeak}. By assumption all these maps are surjective. This shows that according to Definition~\ref{def:AdmissibleGras} that $\sigma\in \delta^S_{\mathfrak{q}}$ is admissible if and only if it is non-trivial. By Corollary~\ref{cor:GrasMunnierRayClassGroup} and the subsequent discussion, we only have to verify for each $\fp\neq \mathfrak{q}_0\in S_p$
    \begin{align}
    \label{eq:SumIsDirect}
        \delta^S_{\mathfrak{q}_0}\cap \langle{\delta^S_\mathfrak{q}:\mathfrak{q}\neq \mathfrak{q}_0}\rangle =\{1\}.
    \end{align}
    Now, the subgroup $\langle{\delta^S_\mathfrak{q}:\mathfrak{q}\neq \mathfrak{q}_0}\rangle$ corresponds by Kummer duality to the field $K'(\sqrt[p]{V^c_{\smash{\mathfrak{q}}}})$, where we set $V_{\mathfrak{q}_0}^c\coloneq\bigcap_{\mathfrak{q}\neq \mathfrak{q}_0}V_\mathfrak{q}$. Thus, $V_{\mathfrak{q}_0}V_{\mathfrak{q}_0}^c=E^{\mathfrak{p}}$ would imply (\ref{eq:SumIsDirect}). This follows from the surjectivity of $\tau^{\mathfrak{p}}$: 
    
    Indeed, let $\varepsilon\in E^{\mathfrak{p}}$ be arbitrary and write $\tau^{\mathfrak{p}}(\varepsilon)=\sum_{\mathfrak{q}}c_\mathfrak{q}$. Now we choose using the surjectivity of $\tau^{\mathfrak{p}}$ two elements $\alpha,\beta\in E^{\mathfrak{p}}$ such that 
    \begin{align*}
        \tau^{\mathfrak{p}}(\alpha)=c_{\mathfrak{q}_0}\quad \text{and}\quad \tau^{\mathfrak{p}}(\beta)=\tau^{\mathfrak{p}}(\varepsilon)-\tau^{\mathfrak{p}}(\alpha).
    \end{align*}
    Then we have $\tau^{\mathfrak{p}}(\varepsilon)=\tau^{\mathfrak{p}}(\alpha\beta)$. Thus, $\varepsilon$ and $\alpha\beta$ differ by an element $\gamma$ of the kernel of $\tau^{\mathfrak{p}}$, which is contained in $V_{\mathfrak{q}_0}$. Therefore, $\varepsilon=(\alpha\gamma)\beta\in V_{\mathfrak{q}_0}V_{\mathfrak{q}_0}^c$.
\end{proof}
For odd $p$ one has $U_{\mathfrak{p}}/U_{\mathfrak{p}}^p \cong \FF_p^{n_\mathfrak{p}+\theta_\mathfrak{p}}$,
where $n_\mathfrak{p}=[K_\mathfrak{p}:\QQ_p]$ and $\theta_\mathfrak{p}=1$ if $K_\mathfrak{p}$ contains a primitive $p^{\text{th}}$-root of unity and $0$ otherwise. Hence it is possible to estimate the dimension of $U_\fp/U_\fp^p$ as $\FF_p$-vector space. If $p$ is unramified in $K/\QQ$, then it is just the inertia degree $f_\mathfrak{p}$.

If $K/\QQ$ is Galois, then $G=\Gal(K/\QQ)$ acts transitively on $S_p$ and for $\sigma\in G$ we get $\sigma\circ \tau^{\mathfrak{p}}=\tau^{\sigma(\mathfrak{p})}\circ \sigma$. Therefore, $\tau^{\mathfrak{p}}$ is surjective if and only if $\tau^{\sigma(\mathfrak{p})}$ is surjective and moreover, $\tau^{\mathfrak{p}}$ is equivariant under the action of the decomposition group $D(\mathfrak{p}|p)$.
        
The next theorem is a direct consequence of Proposition~\ref{prop:WildResidues} together with Theorem~\ref{thm:CriterionH90RayClassGroup} and Proposition~\ref{prop:stabilizationOfC1}.
\begin{repthm}{mainthm:TauCriterion}
    If the map $\tau^\mathfrak{p}$ is surjective for some place $\mathfrak{p}$ of $K$ above $p$, then each $\ZZ_p$-extension in $\mathcal{F}(K,p)$ has the finite coinvariants property. If additionally to the surjectivity of $\tau^{\mathfrak{p}}$ one has $h_p(K)=1$, then each $\ZZ_p$-extension $K_\infty/K$ with $m_0(K_\infty/K)=0$ has the FCP.
\end{repthm}
\begin{rems}
    \label{rem:TauAndDimensions}
    \begin{enumerate}
        \item The proof of Theorem~\ref{mainthm:TauCriterion} does not necessarily require the result of Proposition~\ref{prop:stabilizationOfC1}. There is a more direct argument, as by Proposition~\ref{prop:WildResidues} any cyclic extension $L/K$ where each prime above $p$ is either totally ramified or totally split has the $\bcl^{S_p}$-Hilbert~90 property.
        \item \label{it:DimensionComparison} Assume that $p$ is odd and unramified in $K/\QQ$. Denote by $(r_1,r_2)$ the signature and by $n$ the degree of $K$. Then by the Dirichlet unit theorem, we have $E^{\mathfrak{p}} \otimes \FF_p\cong \FF_p^{r_1+r_2}$.

        Thus, for any $\mathfrak{p}\in S_p$, the map $\tau^{\mathfrak{p}}$ factors as a linear map $\FF_p^{r_1+r_2}\to \FF_p^{n-f_{\mathfrak{p}}}$. Since $r_1+2r_2=n$, we have $r_1+r_2\geq n-f_{\mathfrak{p}}$ if and only if $r_2\leq f_{\mathfrak{p}}$. This is a necessary condition for the surjectivity of $\tau^{\mathfrak{p}}$. 
        \item The conditions of Theorem~\ref{mainthm:TauCriterion} are not hard to verify using a computer. If $p$ is odd and unramified in $K/\QQ$ we have $U_{\mathfrak{q}}^{(2)}\subseteq U_{\mathfrak{q}}^p$  for any $\mathfrak{q}\mid p$ and hence
        \begin{align*}
            U_{\mathfrak{q}}/ U_{\mathfrak{q}}^p\cong (\mathcal{O}_K/\mathfrak{q}^2)^\times\otimes \FF_p.
        \end{align*}
        The map $\tau_\mathfrak{q}^{\mathfrak{p}}:E^{\mathfrak{p}}\to U_{\mathfrak{q}}/ U_{\mathfrak{q}}^p$ is given by reducing an element of $E^{\mathfrak{p}}$ modulo $\mathfrak{q}^2$ and determining its class in $(\mathcal{O}_K/\mathfrak{q}^2)^\times\otimes\FF_p$.
    \end{enumerate}
\end{rems}
{ Similar to the situation and the proof in \ref{prop:Hilbert 90 and logarithmic class group} the statement for the cyclotomic $\ZZ_p$-extension can be shown more directly using $p$-adic class field theory. The following statement and the idea for the proof were communicated to the author by Jaulent.
\begin{prop}
\label{prop:tau for cyclotomic and logarithmic class group}
    Let $K$ be a number field and $p$ a prime such that $\mathcal{F}(K,p)$ contains $K_\infty^{\rm cy}$, e.g., if $p$ is unramified in $K/\QQ$. If $\tau^{\mathfrak{p}}$ is surjective, then the Gross--Kuz'min Conjecture holds for $K$ and $p$.
\end{prop}
\begin{proof}
    We use the notation of \cite{Jaulent1994}. If $\tau^\mathfrak{p}$ is surjective, then by the Nakayama lemma $\mathcal{E}^{\mathfrak{p}}=\ZZ_p\otimes E^{\mathfrak{p}}$ maps surjectively onto $\bigoplus_{\mathfrak{q}\neq \mathfrak{p}}\mathcal{U}_\mathfrak{q}$. By the exact sequence $\smash{0\to \mathcal{U}_\mathfrak{q}\to \mathcal{R}_\mathfrak{p}\to \ZZ_p\to 0}$, where the last map is given by the standard valuation, and a variant of Dirichlet's unit theorem for logarithmic $S$-units, we see that the map $\mathcal{E}^{S_p}\to \bigoplus_{\mathfrak{q}\neq \mathfrak{p}}\mathcal{R}_\mathfrak{p}$ has finite cokernel. Now let $\mathcal{R}_\mathfrak{q}^*$ be the logarithmic local units of $K_\mathfrak{q}$. We see
    \begin{align*}
        \oplus_{\mathfrak{q}\neq \mathfrak{p}} \mathcal{R}_{\mathfrak{q}} \twoheadrightarrow  \oplus_{\mathfrak{q}\neq \mathfrak{p}}\mathcal{R}_\mathfrak{q}/\mathcal{R}_\mathfrak{q}^*\twoheadrightarrow \stilde{\oplus}_{\mathfrak{q}\mid p}\ZZ_p\mathfrak{p}
    \end{align*}
    and conclude from the exact sequence in the proof of Scolie 3.3 in \cite{Jaulent1994} that the kernel of $\logCl_K\to \cl^{S_p}(K)$ is finite. Therefore $\logCl_K$ is itself finite.
\end{proof}
}

Based on the observation in Remark~\ref{rem:TauAndDimensions} (\ref{it:DimensionComparison}) we see that if $K$ is complex, then almost no prime will satisfy the conditions of Theorem~\ref{mainthm:TauCriterion}. Nevertheless, if the primes in $K^+$ above $p$ are inert in $K/K^+$ then Proposition~\ref{prop:H90Transfer} is applicable. Hence in this situation, if the conditions of Proposition~\ref{mainthm:TauCriterion} are satisfied for $K^+$, then it follows that for any $\ZZ_p$-extension $K_\infty^+/K^+$ in $\mathcal{F}(K,p)$, also $KK_\infty^+/K$ has the FCP. Note that assuming the Leopoldt Conjecture, this technique is only applicable to the cyclotomic $\ZZ_p$-extension, as this the only $\ZZ_p$-extension of $K^+$.
\subsection{The imaginary quadratic case}
\label{sec:ImaginaryQuadraticFermat}
Assume now that $K$ is imaginary quadratic and $p$ is unramified in $K/\QQ$. If $p$ is inert, then the conditions of Theorem~\ref{mainthm:TauCriterion} are clearly satisfied and whether a $\ZZ_p$-extension is contained in $\mathcal{F}(K,p)$ depends on the intersection with the $p$-Hilbert class field $H_K(p)$. For example, if $\cl(K)_{(p)}=0$, then every $\ZZ_p$-extension is contained in $\mathcal{F}(K,p)$.

Thus from now on we assume that $p$ is split and $\mathfrak{p}_1$ and $\mathfrak{p}_2$ are the primes above $p$ in $K$. A generator of the torsion-free part of $E^{\mathfrak{p}_1}$ is a generator of $\mathfrak{p}_1^r$, where $r$ is the order of $\mathfrak{p}_1$ in $\cl(K)$. Denote such a generator by $x$. Thus $\tau^{\mathfrak{p}_1}$ is surjective if and only if $x^{p-1}\not\equiv 1\pmod{\mathfrak{p}_2^2}$. This quantity is related to Fermat quotients --- these are numbers of the form $(a^{p-1}-1)/p\in \ZZ$ where $a\in \ZZ$. Here we only restrict ourselves to the valuation of those Fermat quotients, in particular we are interested in $\delta(x)\coloneq v_{\mathfrak{p}_2}(x^{p-1}-1)-1$. Hence $\tau^{\mathfrak{p}_1}$ is surjective if and only if $\delta(x)=0$.

In the cyclotomic $\ZZ_p$-extension the FCP is equivalent to the finiteness of the logarithmic class group. In Theorem 2.2 of \cite{Gras2024} Gras gives an explicit method to compute the size of $\smash{\logCl_K}$ in terms of $\delta(x)$. In particular, if we assume that $\cl(K)_{(p)}=0$, then $|\smash{\logCl_K}|=p^{\delta(x)}$. Thus in some sense our theorem can be seen as a type of higher-dimensional variant of this theorem without the quantitative aspect. 
\begin{quest}
    Is it possible to make the statement of Theorem~\ref{mainthm:TauCriterion} quantitative to compute the size of $\logCl_K$?
\end{quest}
{ For the anticyclotomic $\ZZ_p$-extension $K^{\rm ac}_\infty/K$ we note that $K^{\rm ac}_\infty\in \mathcal{F}(K,p)$ if and only $K^{\rm ac}_\infty\cap H^{\smash{S_p}}_K=K$, which happens for example when $\cl_{[p]}(K)=0$, which is true for all but finitely many primes $p$. Thus we note that if $\tau^{\mathfrak{p}}$ is surjective and $\cl_{[p]}(K)=0$, then $K_\infty^{\rm cy}/K$ and $K_\infty^{\rm ac}/K$ satisfy the FCP simultaneously.

This special case also shows that the problem concerning the surjectivity of $\tau^\mathfrak{p}$ is in the spirit of the Wieferich prime problem: If $\mathfrak{p}_1$ does not satisfy the Wieferich property with respect to $x$, then the conditions of Theorem~\ref{mainthm:TauCriterion} are satisfied. In Section~\ref{sec:Probability} we present a heuristic view inspired by Wieferich prime problem. It predicts that these primes have (natural) density $1$.

Assume that $p\nmid h_K$ and $p$ is split in $K$ as $(p)=\mathfrak{p}_1\mathfrak{p}_2$, then one can always construct a $\ZZ_p$-extension $K_\infty/K$ such that $\mathfrak{p}_2$ is completely split in $K_\infty/K$ and $\mathfrak{p}_1$ is totally ramified in $K_\infty/K$. In particular $K_\infty\in \mathcal{F}(K,p)$. This extension therefore has the finite coinvariants property if $\tau^{\mathfrak{p}}$ is surjective.
}

\subsection{The real quadratic case and \texorpdfstring{$p$}{}-rationality}
For a prime $p$, we say that $K$ is $p$-\emph{rational} if $G_p\coloneq\Gal(K_p/K)$ is a free pro-$p$ group, where $K_p$ is the maximal pro-$p$ extension of $K$ unramified outside the places above $p$ (cf. \cite{Movahhedi1990}).

The $p$-rationality of $K$ is equivalent to the validity of the following two conditions:
\begin{enumerate}
    \item The Leopoldt Conjecture holds for $K$ and $p$;
    \item The $G_p^{ab}$ is torsion-free, which is equivalent to $\rmH^1(G_p,\QQ_p/\ZZ_p)^\vee$ being torsion-free.
\end{enumerate}
We now restrict ourselves to the case where $K$ is real quadratic. Then the Leopoldt Conjecture is satisfied for $K$ and every $p$. Thus, the $p$-rationality of $K$ is equivalent to the vanishing of the torsion subgroup $\mathscr{T}_p$ of $\rmH^1(G_p,\ZZ_p)^*$. It is finite by the Leopoldt Conjecture and thus the vanishing is equivalent to $v_p(|\mathscr{T}_p|)=0$.

We denote by $\tau$ the non-trivial automorphism of $K$ and let $\theta:\{1,\tau\}\to \CC^\times$ be the unique non-trivial complex character of $\Gal(K/\QQ)$. Let $\varepsilon$ be the fundamental unit of $\mathcal{O}_K$. Then we have $v_p(|\mathscr{T}_p|)=v_p(\mathrm{Reg}^\theta_p(\varepsilon))$ by \cite[Sec. 8.6]{Gras2016} for $p$ sufficiently large (it suffices that $p$ does not divide the class number of $K$). 

Let us further assume that $p$ is an odd prime split in $K$ and write $(p)=\mathfrak{p}_1\mathfrak{p}_2$. Then we can compute 
\begin{align*}
    \mathrm{Reg}^\theta_p(\varepsilon)=\frac{-1}{p}\log_p(\varepsilon^2)\equiv \frac{-2}{p-1} \alpha_p(\varepsilon)\pmod{p}
\end{align*}
by \cite[Sect. 2.2.2 and 2.1.1]{Gras2016}, where $\alpha_p(\eta)$ the \emph{generalized Fermat quotient} defined by $\eta^{p-1}=1+p\alpha_p(\eta)$. Note that the definition would differ slightly if $p$ were not split in $K/\QQ$. Now $\mathrm{Reg}^\theta_p(\varepsilon)$ has $p$-valuation $0$ if and only if $\varepsilon^{p-1}\not\equiv 1\pmod {p^2}$, which is equivalent to $\varepsilon_{\mathfrak{p}_1}\not \in U^{p}_{\mathfrak{p}_1}$. This implies that $\tau^{\mathfrak{p}_2}$ is surjective. Therefore, if $K$ is $p$-rational, then it also satisfies the conditions of Theorem~\ref{mainthm:TauCriterion}. The other implication is not true, as the example $K=\QQ(\smash{\sqrt{37}})$ with $p=7$ shows. Here $\varepsilon^6-1\equiv 0\pmod{49}$, but $\tau^{\mathfrak{p}}$ surjective for $\mathfrak{p}$ any prime ideal above $7$.

We conducted experiments for other number fields and there were examples of primes $p$ for which the field was $p$-rational but $\tau^{\mathfrak{p}}$ not surjective and vice versa.

\section{Some statistics and numerical experiments regarding the conditions of Theorem~\ref{mainthm:TauCriterion}}
\label{sec:numerics and statistics}
{
The criterion established in the previous section shows that it is sufficient to check the surjectivity of $\tau^{\fp}$ in order to conclude the $\bcl^{S_p}$-Hilbert~90 property for the first step of several $\ZZ_p$-extensions at once.

There is no obvious structure controlling properties of $\tau^{\fp}$ as a linear map of $\FF_p$-vector spaces of the appropriate dimensions, and in numerical experiments we found that $\tau^{\fp}$ seems to behave \lq\lq randomly\rq\rq\ in most of the cases. We show in Section~\ref{sec:ExceptionalCases} that for CM fields, there are obstructions coming from the complex conjugation.

We start by recalling some basics from the theory of random matrices over finite fields and then use these results to propose a heuristic for the \lq\lq probability\rq\rq\ that $\tau^{\fp}$ is surjective. This leads in special cases to the heuristic usually applied to Wieferich primes. Finally we will present numerical evidence for this heuristic in Section~\ref{sec:numerics}.
}
\label{sec:Probability}
\subsection{Statistical Methods on ranks of random matrices}
\label{sec:StatisticsOnRank}
We first review some techniques from random matrix theory over finite fields and apply them afterwards to our problem. Let $k$ and $n$ be integers and $p$ a prime, then we denote by $\mathbf{T}_p(n,k)\sim \mathrm{Unif}(\FF_p^{n\times k})$ the uniform distribution of matrices in $\FF_p^{n\times k}$, that is, for any subset $U\subseteq  \FF_p^{n\times k}$, we have $\Pr(\mathbf{T}_p(n,k)\in U)=|U|/p^{n\cdot k}$.

Let $\mathbf{R}_p(n,k)$ be the random variable associated to $\FF_p^{n\times k}\to \NN_0$, $A\mapsto \rk(A)$. Clearly $\Pr(\mathbf{R}_p(n,k)<n)=1$ if $k<n$ and by a counting argument we see that for $k\geq n$
\begin{align}
\label{eq:probability for rank}
    \Pr(\mathbf{R}_p(n,k)<n)=1-\prod_{i=0}^{n-1} \left(1-p^{i-k}\right)\leq 1- \left(1-p^{n-1-k}\right)^n.
\end{align}
For a positive real number $x$, we introduce another random variable $\mathbf{P}_{n,k}(x)$ on the set of sequences $(A_p)_p$ with $A_p\in \FF_p^{n\times k}$ indexed by primes $p\leq x$ with discrete uniform distribution, which sends the sequence $(A_p)_p$ to $|\{p\leq x:\rk(A_p)<n\}|$.

Its expected value is
\begin{align*}
    \mathbb{E}[\mathbf{P}_{n,k}(x)]=\sum_{p\leq x}\Pr(\mathbf{R}_p(n,k)<n).
\end{align*}
One can estimate the asymptotic behavior of $\mathbb{E}[\mathbf{P}_{n,k}(x)]$ in terms of $n$ and $k$ using \eqref{eq:probability for rank} as follows.
\begin{prop}
\label{prop:BoundForExpectedValue}
    Let $k$ and $n$ be integers such that $k>n$. Then
    \begin{align*}
        \mathbb{E}[\mathbf{P}_{n,k}(x)]=O(1)\quad \text{and} \quad  
        \mathbb{E}[\mathbf{P}_{n,n}(x)]=O( \log\log x).
    \end{align*}
    The implicit constants depend only on $n$.
\end{prop}
\begin{heu}
\label{heu:ProbilityTauSurjective}
    Let $K$ be a number field of signature $(r_1,r_2)$ and degree $n$, and $p$ be an odd prime, which is unramified in $K/\QQ$. Let $\mathfrak{p}$ be a prime above $p$, then the probability that $\tau^\mathfrak{p}$ is surjective equals 
    \begin{align*}
        \Pr(\mathbf{R}_p(n-f,r_1+r_2)=n-f)
    \end{align*}
    where $f=f(\mathfrak{p}|p)$ is the inertia degree of $\mathfrak{p}$.
\end{heu}
Heuristic~\ref{heu:ProbilityTauSurjective} is clearly correct if $r_2>f(\mathfrak{p}|p)$ for trivial reasons by Remark~\ref{rem:TauAndDimensions} (\ref{it:DimensionComparison}), but it does not hold in all cases. In most of them, the reason was the existence of a certain group action (see Section~\ref{sec:ExceptionalCases}). We believe that this should be the only reason for Heuristic~\ref{heu:ProbilityTauSurjective} to fail. The only counterexamples we were able to find using numerical computations and theoretical arguments are totally imaginary.
\begin{thm}
\label{thm:ProbabilisticProofOfGross}
    Assume that Heuristic~\ref{heu:ProbilityTauSurjective} is valid for the totally real field $K$, then with probability $1$ there are at most finitely many primes $p$ such that the conditions of Theorem~\ref{mainthm:TauCriterion} are not satisfied for the pair $(K,p)$.
\end{thm}
\begin{proof}
    When we write $p$, we mean an odd prime unramified in $K/\QQ$, Furthermore, $\mathfrak{p}$ is an arbitrary prime of $K$ above $p$. By Proposition~\ref{prop:BoundForExpectedValue} we have for $f=f(\mathfrak{p}|p)$ 
    \begin{align*}
        \sum_{p\leq x}\Pr(\rk(\tau^{\mathfrak{p}})<n-f)&\overset{\ref{heu:ProbilityTauSurjective}}\leq \sum_{p\leq x}\Pr(\mathbf{R}_p(n-1,n)<n-1))=\mathbb{E}[\mathbf{P}_{n-1,n}(x)]=O(1).
    \end{align*}
    So by the Borel--Cantelli lemma, with probability $1$ there are only finitely many $p$ such that $\tau^{\mathfrak{p}}$ is not surjective. This finishes the proof.
\end{proof}
Since our numerical computations, which we will present in Section~\ref{sec:StatisticsOnRank}, seem to align quite well with Heuristic~\ref{heu:ProbilityTauSurjective}, we propose the following conjecture:
\begin{conj}
\label{conj:FailiureBound}
    Let $K$ be a totally real number field. There are at most finitely many primes $p$, unramified in $K/\QQ$ for which the map $\tau^{\mathfrak{p}}$ is not surjective for all $\mathfrak{p}$ above $p$. 
\end{conj}
This conjecture extends and implies a conjecture by Jaulent \cite[Conj. 11]{Jaulent2017}, which claims that for a totally real quadratic field, the number of primes $p$, for which $\logCl_K$ is non-trivial, is finite. In fact, if $p$ does not divide $h_K$ and $K_1/K$ has the $\bcl_{[p]}$-Hilbert~90 property, then by the second statement of Proposition~\ref{prop:stabilizationOfC1} we have $\logCl_K=1$. The more general statement, that this holds for each totally real field, is sometimes seen as a \lq\lq folklore conjecture\rq\rq, but to our best knowledge it is made nowhere explicit in the literature.
        
In the same paper Jaulent computes the smallest primes $p$ for which the logarithmic class group $\logCl_K$ is non-trivial for certain quadratic imaginary fields (see Example 9). These agree with the smallest primes $p$, such that $\tau^{\mathfrak{p}}$ is not surjective, except in the case $K=\QQ(\sqrt{-31})$ and $p=2$, as here the contribution comes from the class group of $K$.
\subsection{Some exceptional cases}
\label{sec:ExceptionalCases}
Heuristic~\ref{heu:ProbilityTauSurjective} fails for some classes of number fields due to structural reasons. All the number fields that we found by computer search are totally imaginary. The next proposition gives a criterion for when the conditions of Theorem~\ref{mainthm:TauCriterion} cannot be satisfied for CM fields. 
\begin{prop}
\label{prop:TauFailCM}
    Let $K$ be a CM field and $L$ the Galois closure of $K/\QQ$. We denote by $G=\Gal(L/\QQ)$ and $H=\Gal(L/K)$. The field $L$ is also CM and let $\gamma$ be the complex conjugation in $L$. Let $p$ be an odd prime number which is unramified in $K/\QQ$ such that there are at least two places above $p$ in $K$. Let $\mathfrak{P}$ be a prime of $L$ above $p$. If $\gamma\in \langle H,D(\mathfrak{P}|p)\rangle$, then the map $\tau^{\mathfrak{p}}$ is not surjective, where $\mathfrak{p}$ is the place of $K$ below $\mathfrak{P}$.
\end{prop}
\begin{proof}
    Let $\mathfrak{p}_1$,...,$\mathfrak{p}_g$ be the primes of $K$ above $p$. After reordering the places $\mathfrak{p}_i$, we can assume that $\mathfrak{p}=\mathfrak{p}_1$. We set $K_0\coloneq K\cap L^{D(\mathfrak{P}|p)}$, then we have $K_0\subseteq K^+$, as the subgroup $\Gal(L/K_0)=\langle H,D(\mathfrak{P}|p)\rangle$ contains $\gamma$, which restricted to $K$ is the complex conjugation of $K$. This implies that $\mathfrak{p}_1$ is inert in the extension $K/K_0$ and therefore also in $K/K^+$. 
     
    The subgroup $C_2\cong \langle \gamma\rangle$ of $G$ acts on $E_K^{\mathfrak{p}_1}$, as $\gamma$ fixes $\mathfrak{p}_1$. This action descends to $A\coloneq E_K^{\mathfrak{p}_1}\otimes \FF_p$. Since the group algebra $\FF_p[C_2]$ is semi-simple, $A$ decomposes as $A^+\oplus A^-$, where $A^+$ is the maximal submodule of $A$, on which $\gamma$ acts as $1$ and $A^-$ the submodule of $A$, on which $\gamma$ acts as $-1$. 
    
    Next, we are going to show that $A^-=0$. Consider the inclusion map $E^{\mathfrak{p}_1}_{K^+}\to E^{\mathfrak{p}_1}_K$. Since $\mathfrak{p}_1$ is inert in $K/K^+$ we can use argument as in \cite[Thm. 4.12]{Washington1997} to show that the order of the cokernel of $E^{\mathfrak{p}_1}_{K^+}\to E^{\mathfrak{p}_1}_K$ divides $|\mu(K)|$. Since $p$ is unramified in $K/\QQ$ it is coprime to $2\cdot |\mu(K)|$. Thus, the induced map 
    \begin{align*}
        E^{\mathfrak{p}_1}_{K^+}\otimes \FF_p\to E^{\mathfrak{p}_1}_K \otimes \FF_p
    \end{align*}
    is an isomorphism. Note that $\gamma$ acts trivially on the left side and the inclusion commutes with the $\gamma$ action. This implies that $\gamma$ acts trivially on $E_K^{\mathfrak{p}_1}\otimes \FF_p$ and $A^-=0$. Since $\tau^{\mathfrak{p}_1}$ commutes with the action of $\gamma$, we are done if $\gamma$ acts non-trivially on $\bigoplus_{i=2}^g U_{\mathfrak{p}_i,K}/U_{\mathfrak{p}_i,K}^p$, which is the codomain of $\tau^{\mathfrak{p}_1}$. If $\gamma$ permutes the places $\mathfrak{p}_i$ non-trivially, we are done. So we can assume $\gamma(\mathfrak{p}_i)=\mathfrak{p}_i$ for all $i$. But this implies that all $\mathfrak{p}_i$ are inert in $K/K^+$.

    We now show that $\gamma$ acts non-trivially on $M\coloneq U_{\mathfrak{p}_2,K}/U_{\mathfrak{p}_2,K}^p$. This will complete the proof. We compute the invariants $M_{\gamma}$ by the following exact sequence:
    \begin{equation*}
        \begin{tikzcd}[column sep=small]
            0\arrow[r]&(U_{\mathfrak{p}_2,K}^p)_\gamma\arrow[r]&(U_{\mathfrak{p}_2,K})_\gamma\arrow[r]&M_\gamma  \arrow[r]&0
        \end{tikzcd}
    \end{equation*}
    The exactness follows, as the connecting homomorphism $M_\gamma\to \rmH^1(C_2,U_{\mathfrak{p}_2,K}^p)$ maps a $p$-torsion abelian group into a $2$-torsion abelian group and is therefore $0$. This shows that $M_\gamma$ is a quotient of $(U_{\mathfrak{p}_2,K})_\gamma\cong U_{\mathfrak{p}_2,K^+}$. As $U_{\mathfrak{p}_2,K^+}^p\subseteq (U_{\mathfrak{p}_2,K}^p)_\gamma$, we get that $U_{\mathfrak{p}_2,K^+}/ U_{\mathfrak{p}_2,K^+}^p \rightarrow M_\gamma$
    is a surjection. But the left hand side has $\FF_p$ dimension $f_2/2$. Since $M$ has $\FF_p$-dimension $f_2$, $M_\gamma\neq M$ and $\gamma$ acts non-trivially on $M$.
\end{proof}
If $K/\QQ$ is Galois, then Proposition~\ref{prop:TauFailCM} reads as follows.
\begin{cor}
\label{cor:TauFailCMGalois}
    Let $K$ be a CM-field with complex conjugation $\gamma$ such that $K/\QQ$ is Galois. Let $p$ be an odd prime number, which is unramified in $K/\QQ$ and such that there are at least two places above $p$ in $K$. Let $\mathfrak{p}$ be a prime of $K$ above $p$. If $\gamma\in D(\mathfrak{p}|p)$, then the map $\tau^\mathfrak{p}$ is not surjective.
\end{cor}
In combination with the Chebotarev Density theorem, Corollary~\ref{cor:TauFailCMGalois} allows to make more precise predictions about the number of primes $p$ satisfying the conditions of Theorem~\ref{mainthm:TauCriterion}. We demonstrate this with an example where $\Gal(K/\QQ)\cong C_6\times C_2$.
\begin{exmp}
\label{exmp:exceptionalCaseForC2C6}
    Let $K$ be a CM-field with $G=\Gal(K/\QQ)\cong C_6\times C_2$ and $\gamma\in \Gal(K/\QQ)$ the complex conjugation. Let $p$ be an odd prime  unramified in $K/\QQ$. For dimension reasons if $p$ is a prime such that there are more than $2$ primes above $p$, the map $\tau_p$ cannot be surjective (see Remark~\ref{rem:TauAndDimensions}). So assume that $(p)=\mathfrak{p}_1\mathfrak{p}_2$. Then $D(\mathfrak{p}_i)\cong C_6$. Since $G$ contains $3$ cyclic subgroups of order $6$ and for any cyclic subgroup $U$ of order $2$, there are exactly two cyclic subgroups of order $6$ not containing $U$.

    As $D(\mathfrak{p}_i)$ is generated by the Frobenius of $p$, there are $2\varphi(6)=4$ possibilities for the Frobenius of $p$, for which the conditions of Corollary~\ref{cor:TauFailCMGalois} are not satisfied. By the Chebotarev density theorem and Heuristic~\ref{heu:ProbilityTauSurjective}, we would predict that $4/12=1/3$ of the primes that satisfy the conditions of Theorem~\ref{mainthm:TauCriterion}. Whereas previously we would estimate that roughly $50\%$ of primes satisfy the conditions. Table~\ref{tab:improved estimates for CM} shows for different number fields $K=\QQ[x]/(f)$ estimates and the real of the number of primes up to $10^5$, for which $\tau_p$ is surjective.
    \begin{table}[h!]
        \centering
        \begin{tabular}{c|c|c|c|c}
            Polynomial & Previous est. & New est. & Total primes & Satisfied cases \\
            \hline 
            \hline
            $x^{12}-x^6+1$ & $4820$ &$3213$ & $9590$ &$3207$ \\
            $x^{12}-8x^6+64$ & $4799$ & $3199$  & $9590$ & $3196$\\
            $x^{12}+13x^8+26x^4+1$& $4802$ & $3201$ & $9590$ &  $3208$
        \end{tabular}
        \caption{New estimates based on Proposition~\ref{prop:TauFailCM}}
        \label{tab:improved estimates for CM}
    \end{table}
    
\end{exmp}

\begin{rem}
    There seem to be further exceptional cases of number fields $K$, which do not follow Heuristic~\ref{heu:ProbilityTauSurjective} but neither satisfy the condition of Proposition~\ref{prop:TauFailCM}. One of these examples is $K=\QQ[x]/(f)$, where 
    \begin{align*}
        f=x^8-3x^7+3x^6+x^5-2x^4-3x^3+7x^2-4x+1.
    \end{align*}
    The group $\Gal(f/\QQ)$ is a semidirect product $Q_8\rtimes_\varphi C_2$, where $\varphi$ is a non-trivial homomorphism, $C_2\to \Aut(Q_8)\cong S_4$ such that the composition with $\Aut(Q_8)\to \mathrm{Out}(Q_8)\cong S_3$ is trivial. One can show that any such choice gives isomorphic groups.
    
    The number of primes $\leq 10^5$ that satisfy the conditions of Theorem~\ref{mainthm:TauCriterion} is $3608$, while Heuristic~\ref{heu:ProbilityTauSurjective} predicts the number to be $4815$ out of $9589$. Neither $K$ nor its Galois closure is CM, but totally imaginary.
\end{rem}

\subsection{Numerical experiments related to Heuristic~\ref{heu:ProbilityTauSurjective}}
\label{sec:numerics}
In this section, we present the results of numerical experiments studying Heuristic~\ref{heu:ProbilityTauSurjective}. In the first part, we fix the field $K$ and let $p$ vary. The advantage of this approach is that it is natural and gives (in the totally real case) an intuition for the conjectural statement of Theorem~\ref{thm:ProbabilisticProofOfGross}. There are clear  disadvantages of this strategy. The first one is, that we are only able to check consequences of the heuristic but not the statement itself. Secondly, since the number of cases, where the condition fails is very small we cannot draw conclusions about the precise asymptotic behavior.

In the second part, we fix the prime $p$ and vary the field $K$ in a family. The advantages and disadvantages are opposite to the previous case, as here we can draw conclusions about the precise value of the claimed probability, but there is no canonical way to construct families of number fields. We choose an approach based on Hilbert's irreducibility theorem, but there are various other methods to group number fields into families. 

The data presented here was generated using the computer algebra system OSCAR. The source code can be found under \cite{Feuerpfeil2026Code}.

\subsubsection{Fixing \texorpdfstring{$K$}{} and varying \texorpdfstring{$p$}{}}
\label{sec:FixingKAndVaryp}
Fix a number field $K$ of degree $n$ and signature $(r_1,r_2)$. For abbreviation, we set $r\coloneq r_1+r_2$. Then for a nonnegative integer $k$ and $x>0$, we define
\begin{align*}
    R_p(k)\coloneq\prod_{i=0}^{k-1} \left(1-p^{i-r}\right)\qquad \text{and}\qquad
    \mathcal{E}(x)\coloneq\sum_{(*)}R_p(n-\max_{\mathfrak{p}\in S_p}f_{\mathfrak{p}}),
\end{align*}
where $(*)$ indicates the primes $p\leq x$ that are unramified in $K/\QQ$. Thus, $\mathcal{E}(x)$ is the number of primes up to $x$ for which we expect the conditions of Theorem~\ref{mainthm:TauCriterion} to hold. Let $\Pi_{nr}(x)$ be the set of the primes $p\leq x$, which are unramified in $K/\QQ$. Then we set
\begin{align*}
    \pi_{nr}(x)&\coloneq|\Pi_{nr}(x)|\\
    \pi_{nr}^f(x)&\coloneq|\{p\in \Pi_{nr}(x): \max_{\mathfrak{p}\in S_p}f_\mathfrak{p}<r_2\}|\\
    \mathcal{R}(x)&\coloneq|\{p\in \Pi_{nr}(x):p\text{ satisfies Theorem~\ref{mainthm:TauCriterion}}\}|\\
    \overline{\mathcal{R}}(x)&\coloneq\pi_{nr}(x)-\pi_{nr}^f(x)-\mathcal{R}(x).
\end{align*}
Thus $\overline{\mathcal{R}}(x)$ describes the number of primes for which the conditions of Theorem~\ref{mainthm:TauCriterion} do not fail for obvious reasons but do not satisfy them nonetheless. 

We now give some concrete examples of number fields $K$ and their associated data. We give defining polynomials of the number fields.  We present the data organized by $r_2$ and put special emphasis on the cases $r_2=0$ and $r_2=1$. For all the computations, we set $x=10^5$.
\begin{itemize}
    \item[$\mathbf{r_2=0}$:] In this case $\pi_{nr}^f(x)=0$. The results are depicted in Table~\ref{tab:0r2Experiment}.
    \begin{table}[h!]
        \centering
        \begin{tabular}{c|c|c|c|c}
           $f$ & $\pi_{nr}(x)$ &  $\mathcal{R}(x)$ &$ \mathcal{E}(x)$& $\overline{\mathcal{R}}(x)$  \\
           \hline
           \hline
           $t^5-t^4-34t^3+285t+285$&$9587$&$9587$&$9587.0$&$0$\\
          
           $t^5-55t^3+605t-22$&$9587$&$9587$&$9587.0$&$0$\\
          
           $t^6-10t^4+24t^2-8$&$9590$&$9589$&$9590.0$&$1$\\
          
           $t^7-7t^5+13t^3-5t-1$&$9591$&$9591$&$9591.0$&$0$\\
          
           $t^8-13t^6+44t^4-17t^2+1$&$9589$&$9589$&$9588.9$&$0$\\
        \end{tabular} 
        \caption{Numerical experiments for totally real fields}
        \label{tab:0r2Experiment}
    \end{table}  
    
    \item[$\mathbf{r_2=1}$:] Also here we have $\pi_{nr}^f(x)=0$. The data is listed in Table~\ref{tab:1r2Experiment}.
    \begin{table}[h!]
        \centering
        \begin{tabular}{c|c|c|c|c}
           $f$ & $\pi_{nr}(x)$ &  $\mathcal{R}(x)$ &$ \mathcal{E}(x)$& $\overline{\mathcal{R}}(x)$  \\
           \hline
           \hline
            $t^5 - 41(t-1)(t-2)(t-3)$&$9588$&$9587$&$9587.9$&$1$\\
           
            $t^6 - t^5 - t^4 - 2t^3 + 2t^2 + 3t - 1$&$9590$&$9589$&$9589.9$&$1$\\
           
            $t^7-t^5-4t^3+4t-1$&$9591$&$9591$&$9590.9$&$0$\\
           
            $t^8-6t^4+5t^2-1$&$9589$&$9589$&$9588.9$&$0$\\
           
            $t^8-4t^6+t^4+5t^2-2$&$9589$&$9589$&$9588.9$&$0$\\
            
       \end{tabular} 
       \caption{Numerical experiments in the case $r_2=1$}
        \label{tab:1r2Experiment}
    \end{table}
    
   For the case $K=\QQ(\sqrt[4]{11})$ we extended the computations to $x=10^7$. Under the $664577$ unramified primes, we only found three primes, for which the conditions did not apply. These primes are $19$, $7603$ and $1162223$.
   
   Using the criterion from Theorem~\ref{thm:CriterionH90RayClassGroup} one can verify that for all three primes, the first step in the cyclotomic $\ZZ_p$-extension $K_\infty^{\rm cy}/K$ has the $\bcl^S_{(p)}$-Hilbert~90 property.
    \item[$\mathbf{r_2\geq 2}$] Here we collect some examples where $r_2\geq 2$. The computational data can be found in Table~\ref{tab:2r2Experiment}.
    \begin{table}[h!]
        \centering
        \begin{tabular}{>{\centering}p{4cm}|c|c|c|c|c|c}
           $f$ & $r_2$ & $\pi_{nr}(x)$ & $\pi_{nr}^f(x)$ & $\mathcal{R}(x)$ &$ \mathcal{E}(x)$ & $\overline{\mathcal{R}}(x)$  \\
           \hline
           \hline
           $t^6 - (t^2-2)(t^3+t^2)+3t - 1$&$2$&$9589$&$9$&$9580$&$9579.7$&$0$\\
          
           $t^{7} - t^{4} - t^{3} + t^{2} + 1$&$3$&$9592$&$430$&$9162$&$9161.7$&$0$\\
          
            $t^8+t^6-t^3-t^2-1$&$3$&$9589$&$182$&$9407$&$9406.8$&$0$\\
           
            $t^8-t^7-t^3+2t^2-t+1$&$4$&$9589$&$2897$& $6689$&$6691.4$&$3$\\
           
            $t^9-t^7-t^6+t^4-t^3+t+1$&$4$&$9591$&$780$&$8811$&$8810.5$&$0$
        \end{tabular}   
        \caption{Numerical experiments for fields with at least two pairs of complex embeddings}
        \label{tab:2r2Experiment}
    \end{table}
    
\end{itemize}

\subsubsection{Fixing $p$ and varying $K$}
To produce a family of fields, we use Hilbert's
Irreducibility Theorem (\cite{Hilbert1892}). See, for example, \cite[Chap. 9]{Lang1983} for a comprehensive introduction and a proof of the following theorem (in a more general form):
\begin{thm}[Hilbert's Irreducibility Theorem]
    Let $g\in \QQ[X,t]$ be an irreducible polynomial in two variables. Then the polynomial $g(X,\lambda)\in \QQ[X]$ is irreducible for infinitely many $\lambda\in \QQ$. In this case $\Gal(g(X,\lambda)/\QQ)$ is a (transitive) subgroup of $\Gal(g(X,t)/\QQ(t))$.
\end{thm}
We first outline the general strategy: Let $g\in \QQ[X,t]$ be a polynomial in two variables of $X$-degree $n$. We fix two integers $a<b$ and set $\Lambda_{a,b}\coloneq[a,b]\cap \ZZ$. In all the cases, we choose $a$ and $b$ such that  $g(X,\lambda)$ for $\lambda\in \Lambda_{a,b}$ has a constant number of real roots in this set, in particular the signatures of all the fields listed in one table agree. For a fixed prime $p$, we determine the number of $\lambda\in \Lambda_{a,b}$ such that $g(X,\lambda)$ is irreducible and $K_{g,\lambda}\coloneq\QQ[X]/(g(X,\lambda))$ satisfies the conditions of Theorem~\ref{mainthm:TauCriterion}. We compare it to the predicted number by Heuristic~\ref{heu:ProbilityTauSurjective}. 

Define the following sets:
\begin{align*}
    \Lambda^{\text{irr}}_{a,b}&\coloneq\{\lambda\in \Lambda_{a,b}:g(X,\lambda)\text{ is irreducible and }p\text{ unramified in }K_{g,\lambda}/\QQ\}\\
    \Lambda^{p}_{a,b}&\coloneq\{\lambda\in \Lambda_{a,b}^{\text{irr}}:K_{g,\lambda}\text{ statisfies the conditions of Theorem~\ref{mainthm:TauCriterion} with respect to }p\}
\end{align*}
Similar to the treatment in Section~\ref{sec:FixingKAndVaryp} we define define the expected number of valid cases $\mathcal{E}(a,b)$ to be
\begin{align*}
    \mathcal{E}(a,b)\coloneq\sum_{\lambda\in \Lambda_{a,b}^{\text{irr}}}R_p(\deg_X(g)-f_p^{\max}).
\end{align*}
Note that $\deg_X(g)$ is the extension degree of $K_{g,\lambda}$ for $\lambda\in \Lambda^{\rm irr}_{a,b}$.

As polynomials $g$, we use some of the ones computed by Smith in \cite{Smith2000}.
\begin{enumerate}
    \item Let $g(X,t)\coloneq X^3+tX^2+(t-3)X-1$. Then $\Gal(g(X,t)/\QQ(t))\cong C_3$. Therefore we have $\Gal(g(X,\lambda)/\QQ)\cong C_3$ for any $\lambda$ such that $g(X,\lambda)$ is irreducible. Furthermore, in that case $K_{g,\lambda}$ will be totally real. The results of our computations for $a=-1000$ and $b=1000$ can be found in Table~\ref{tab:VaryingFieldC3}.
    \begin{table}[h!]
        \centering
        \begin{tabular}{c|c|c|c}
             $p$& $\#\Lambda_{a,b}^{\text{irr}}$ & $\# \Lambda_{a,b}^{p}$ & $\mathcal{E}(a,b)$\\
             \hline\hline
             $3$&$1408$&$1408$&$1404.54$\\
             $5$&$2001$&$2001$&$1981.9$\\
             $7$&$1439$&$1421$&$1432.3$\\
             $11$&$2001$&$2001$&$1996.1$\\
             $31$&$1873$&$1871$&$1872.4$
        \end{tabular}
        \caption{Statistics for the generic polynomial $X^3+tX^2+(t-3)X-1$}
        \label{tab:VaryingFieldC3}
    \end{table}
    \item Let $g(X,t)\coloneq X^3-3X-t$. Then $\Gal(g(X,t)/\QQ(t))\cong S_3$. For $|\lambda|>2$, the polynomial $g(X,\lambda)$ has exactly one real root. Thus if $g(X,\lambda)$ is irreducible, we will have $\Gal(g(X,\lambda)/\QQ)\cong S_3$. Again for $a=-1000$ and $b=1000$ the results are depicted in Table~\ref{tab:VaryingFieldS3}.
    \begin{table}[h!]
        \centering
        \begin{tabular}{c|c|c|c}
             $p$& $\#\Lambda_{a,b}^{\text{irr}}$ & $\# \Lambda_{a,b}^{p}$ & $\mathcal{E}(a,b)$\\
             \hline\hline
             $2$&$1074$&$1036$&$1015.8$\\
             $5$&$1322$&$1283$&$1290.0$\\
             $7$&$1488$&$1463$&$1464.8$\\
             $37$&$1876$&$1866$&$1868.0$\\
             $101$&$1946$&$1940$&$1942.8$
        \end{tabular}
        \caption{Statistics for the generic polynomial $X^3-3X-t$}
        \label{tab:VaryingFieldS3}
    \end{table}
    \item Let $g(X,t)\coloneq X^4+tX^2+1$. Then $\Gal(g(X,t)/\QQ(t))\cong C_2\times C_2$. For $\lambda < -2$, the polynomial $g(X,\lambda)$ has exactly $4$ real roots, whereas for $\lambda>-2$ it only has complex roots. Since the action of $C_2\times C_2$ on the roots is regular, we will have $\Gal(K_{g,\lambda}/\QQ)\cong C_2\times C_2$ whenever $g(X,\lambda)$ is irreducible. This implies that for $\lambda>-2$ the field $K_{g,\lambda}$ is CM and we have to take into account the exceptional case of Corollary~\ref{cor:TauFailCMGalois}. Since we do not have an analogue of the Chebotarev-Density theorem for this case at hand, it is hard to make good predictions. In Table~\ref{tab:VaryingFieldC2C2} we present the data for $a_1=-2000$ and $b_1=-1$ on the left and for $a_2=1$ and $b_2=2000$ on the right.
    \begin{table}[h!]
        \subfloat{
        \begin{tabular}{c|c|c|c}
             $p$& $\#\Lambda_{a_1,b_1}^{\text{irr}}$ & $\# \Lambda_{a_1,b_1}^{p}$ & $\mathcal{E}(a_1,b_1)$\\
             \hline\hline
             $3$&$933$&$789$&$880.5$\\
             $5$&$1278$&$1199$&$1263.7$\\
             $7$&$1424$&$1350$&$1413.6$\\
             $11$&$1587$&$1559$&$1582.8$\\
             $13$&$1640$&$1618$&$1636.5$
        \end{tabular}}
        \hspace{.3cm}
         \subfloat{\begin{tabular}{c|c|c|c}
             $p$& $\#\Lambda_{a_2,b_2}^{\text{irr}}$ & $\# \Lambda_{a_2,b_2}^{p}$ & $\mathcal{E}(a_2,b_2)$\\
             \hline\hline
             $3$&$998$&$435$&$542.2$\\
             $5$&$1333$&$598$&$973.1$\\
             $7$&$1397$&$676$&$1003.4$\\
             $11$&$1665$&$778$&$1167.5$\\
             $13$&$1713$&$776$&$1281.9$
        \end{tabular}}
        \caption{Statistics for the generic polynomial $X^4+tX^2+1$}
        \label{tab:VaryingFieldC2C2}
    \end{table}
    
    For $\lambda>0$, it seems as if the probability that Theorem~\ref{mainthm:TauCriterion}  is satisfied is roughly $50\%$, whereas the precise prediction is about $75\%$. This fits well with a similar argument, as the one we used in Example~\ref{exmp:exceptionalCaseForC2C6}.
\end{enumerate}

\section*{Acknowledgments}
I would like to thank my Ph.D. advisors, Christian Maire (Université Marie et Louis Pasteur) and Thomas Weigel (Università degli Studi di Milano-Bicocca), for proposing this problem and for their constant support and encouragement throughout this work. I also wish to thank Donghyeok Lim and Ravi Ramakrishna for many stimulating discussions. I am indebted to Georges Gras for his insightful comments and valuable suggestions on an earlier version of this article, which have greatly improved its presentation and style. I would especially like to thank Jean-Fran\c{c}ois Jaulent, for his interest in this work and many insights into $\ell$-adic class field theory.

\noindent I am a member of the Gruppo Nazionale per le Strutture Algebriche, Geometriche e le loro Applicazioni (GNSAGA) of the Istituto Nazionale di Alta Matematica (INdAM).

\noindent I would like to thank an anonymous referee for carefully reading the article and for his/her insightful comments leading to a significant improvement of the article.

\section*{Funding}
During the research, the author received a mobility grant (ANR-17-EURE-0002) by Graduate school EIPHI at Université Marie et Louis Pasteur and Université Bourgogne Europe.

\section*{Data Availability Statement}
All data generated and analyzed during this study were produced using code developed by the author based on the open-source computer algebra system OSCAR (see \cite{OSCAR}). The complete dataset underlying the results is fully presented within this article. The source code used to generate the data is publicly available in the author's Github under \cite{Feuerpfeil2026Code}.

\appendix

\section{Brauer Groups and the \texorpdfstring{$\mathbf{Pic}$}{}-Hilbert~90 Property for Cyclic Covers of Schemes}
\label{app:GeometricProof}
This appendix is logically independent of the rest of the paper and gives an alternative proof of Theorem~\ref{thm:Hilbert90CriterionViaNorms} using more algebro-geometric methods, which applies to a broader context and might have applications in the function-field case. In doing so, we also give another example of a cohomological Mackey functor and generalize some previous results before to cyclic Galois coverings of schemes. 

Fix a finite Galois cover $\pi:Y\to X$ of schemes, with Galois group $G$. Then we define a cohomological Mackey functor $\mathbf{H}^n_{\et}(Y/X)$ by
\begin{equation*}
 \begin{tikzcd}
 \rmH^n_{\et}(Y, \mathbb{G}_{m})\arrow[r, shift left, "N_\pi"]&\arrow[l, shift left, "f_*"]\rmH^n_{\et}(X, \mathbb{G}_{m})
 \end{tikzcd}
\end{equation*}
Here, $N_\pi$ is induced by the norm map $\pi_*\mathcal{O}_Y^*\to \mathcal{O}_X^*$, as defined by Mumford in \cite[Lecture 10]{Mumford1966}. It agrees locally with the map $\alpha\mapsto \prod_{\sigma\in G}\sigma(\alpha)$, which makes it easy to check that it is indeed a cohomological Mackey functor. We abbreviate $\mathbf{H}^1_{\et}(Y/X)$ by $\mathbf{Pic}_{Y/X}$, or, if there is no confusion about the referred cover, by $\mathbf{Pic}$. The notation should not be confused with the Picard scheme.

The Hochschild-Serre spectral sequence in étale cohomology (see \cite[§2 Theorem 2.20]{Milne1980}) now yields the following seven-term sequence:
\begin{equation}
\label{eq:SevenTermSequence}
\begin{split}
 0\to \rmH^1(G, \mathbb{G}_m(Y))&\to \Pic(X)\to\Pic(Y)^G\to \rmH^2(G, \mathbb{G}_m(Y))\to\phantom{ab}\\ \to&\Br(Y/X)\to \rmH^1(G, \Pic(Y))\to \rmH^3(G, \mathbb{G}_m(Y))
\end{split}
\end{equation}
Here we write $\Br(Y/X)$ for the relative cohomological Brauer group:
\begin{align*}
 \Br(Y/X): =\ker(\rmH^2_{\et}(X, \mathbb{G}_m)\overset{f_*}\to \rmH^2_{\et}(Y, \mathbb{G}_m)).
\end{align*}
The seven-term sequence (\ref{eq:SevenTermSequence}) was independently discovered for Galois extensions of rings by Chase, Harrison, and Rosenberg in \cite[Corollary 5.5]{ChaseHarrisonRosenberg1965}, whose methods were based on Amitsur cohomology.

We use this sequence to deduce some properties of the cohomological Mackey functor $\mathbf{Pic}$. One immediately observes that
\begin{align*}
 k^0(G, \mathbf{Pic})&\cong \rmH^1(G, \mathbb{G}_m(Y)) \quad \text{and}\quad  k^1(G, \mathbf{Pic})\cong \ker (\rmH^2(G, \mathbb{G}_m(Y))\to \Br(Y/X)).
\end{align*}
From now on, we assume $G$ to be cyclic. Then we can identify $\smash{\widehat{\rmH}}^i(G,M)$ with $\smash{\widehat{\rmH}}^{i+2n}(G,M)$ for any $n\in \ZZ$. Such an isomorphism depends on the choice of a generator $\sigma$ of $G$.
\begin{prop}
\label{prop:DiagramCommutative}
 Assume that $Y$ and $X$ are locally factorial, regular, and integral. Then the following diagram is (for any choice of generator of $G$) commutative up to sign:
 \begin{equation*}
 \begin{tikzcd}
 \smash{\widehat{\rmH}}^{-1}(G, \Pic(Y))\arrow[r]\arrow[d, "\sim" {rotate=90, anchor=north}]&k^1(G, \mathbf{Pic})\arrow[r, "\sim"]&\rmH^1(G, \mathbb{G}_m(Y))\arrow[d, "\sim" {rotate=90, anchor=north}]\\
 \rmH^1(G, \Pic(Y))\arrow[rr, "d_2^{1, 1}"] &&\rmH^3(G, \mathbb{G}_m(Y)
 \end{tikzcd}
 \end{equation*}
 Here the bottom horizontal arrow is the differential in the Hochschild--Serre spectral sequence and the top left arrow the one induced from the six term sequence (\ref{eq:sixtermsequence}).
\end{prop}
\begin{proof}
 The proof proceeds in exactly the same as the one in \cite[Lemma 4.2]{Beauville2009}. The result of \cite[Proposition 1.1]{Skorobogatov2007} still holds true in this slightly different situation by considering the exact functor $\mathbf{mod}_G\to \mathrm{Sh}(X_{\et})$, $M\mapsto f_*(\underline{\ZZ})\otimes_{\underline{\ZZ[G]}}\underline{M}$, which is left adjoint to $\Gamma(Y,\bl):\mathrm{Sh}(X_{\et})\to \mathbf{mod}_G$.
\end{proof}
González-Avilés noted in \cite[Remarks 6.3 (b)]{Aviles2020} that the statement of Proposition \ref{prop:DiagramCommutative} should be true in a wider generality, but we were unable to find a proof in a more general setting. However, this statement is sufficient for our applications.

Proposition \ref{prop:DiagramCommutative} together with the six term sequence (\ref{eq:sixtermsequence}) for cohomological Mackey functors implies the following theorem:
\begin{thm}
 \label{thm:exactHexagon}
 Let $\pi:Y\to X$ be a Galois cover of locally factorial, regular, and integral schemes with cyclic Galois group $G$. Then the following hexagon is exact:
 \begin{center}
     \begin{tikzpicture}[scale=.8]

\node[xshift=1.6cm] (A) at (60:1.5cm)
    {$\smash{\widehat{\rmH}}^{1}(G,\mathbb{G}_m(Y))$};

\node[xshift=2.0cm] (B) at (0:1.5cm)
    {$c_0(G,\mathbf{Pic})$};

\node[xshift=1.6cm] (C) at (-60:1.5cm)
    {$\smash{\widehat{\rmH}}^{0}(G,\Pic(Y))$};

\node (D) at (-120:1.5cm)
    {$\smash{\widehat{\rmH}}^{0}(G,\mathbb{G}_m(Y))$};

\node[xshift=-.4cm] (E) at (-180:1.5cm)
    {$\Br(Y/X)$};

\node (F) at (-240:1.5cm)
    {$\smash{\widehat{\rmH}}^{1}(G,\Pic(Y))$};

\draw[->] (A) -- (B);
\draw[->] (B) -- (C);
\draw[->] (C) -- (D);
\draw[->] (D) -- (E);
\draw[->] (E) -- (F);
\draw[->] (F) -- (A);

\end{tikzpicture}
 \end{center}
 In particular, we have $c_1(G,\mathbf{Pic})\cong \coker(\smash{\widehat{\rmH}}^{0}(G,\mathbb{G}_m(Y))\to \Br(Y/X))$ and $G$ has $\mathbf{Pic}$-Hilbert~90 property if and only if the map $\smash{\widehat{\rmH}}^{0}(G,\mathbb{G}_m(Y))\to \Br(Y/X)$ is surjective.
\end{thm}
\begin{exmp}
 Let $\mathcal{C}/k$ be a smooth projective curve over an algebraically closed field $k=\overline{k}$ and $\pi:\mathcal{C}'\to \mathcal{C}$ a connected Galois cover of $C$ with cyclic Galois group $G$. Then $\mathcal{C}'/\mathcal{C}$ has the $\mathbf{Pic}$-Hilbert~90 property, as 
 \begin{align*}
 \Br(\mathcal{C}'/\mathcal{C})\hookrightarrow \rmH^2_{\et}(\mathcal{C}, \mathbb{G}_m)\hookrightarrow \rmH^2(k(\mathcal{C}), \mathbb{G}_m)
 \end{align*}
 and the latter one vanishes by Tsen's theorem. In fact, it is sufficient that $\mathrm{cd}_\ell(k)=0$ for all primes $\ell$ dividing $|G|$.

 This applies, for example, to cyclic isogenies of elliptic curves over algebraically closed fields.
\end{exmp}
If $L/K$ is a cyclic extension of number fields and $S$ a set of places containing the ramified and Archimedean ones, then $\pi:\Spec(\mathcal{O}_{L,S})\to \Spec(\mathcal{O}_{K,S})$ is a Galois cover with Galois group $G\coloneq \Gal(L/K)$. A proof of the following lemma can be found in \cite[Appendix Lemma A.1]{Aviles2004}:
\begin{lem}
\label{lem:SesForBrauer}
 Let $\mathcal{O}_{K,S}$ resp. $\mathcal{O}_{L,S}$ be the ring of $S$-integers of $K$ resp. $L$, then there is an exact sequence:
 \begin{equation*}
 \begin{tikzcd}
 0\arrow[r]&\Br(\mathcal{O}_{L, S}/\mathcal{O}_{K, S})\arrow[r]&\bigoplus_{v\in S}\Br(L_v/K_v)\arrow[r, "{\sum \inv_v}"]&\QQ/\ZZ
 \end{tikzcd}
 \end{equation*}
\end{lem}
Now, the proof of Theorem~\ref{thm:Hilbert90CriterionViaNorms} is a simple combination of Theorem~\ref{thm:exactHexagon} with the above lemma and the decomposition of every abelian group involved into their $p$-primary parts.
\bibliography{references}

@article{QuadrelliWeigel2015,
 author = {Quadrelli, C. and Weigel, Th.},
 title = {A group-theoretical version of {Hilbert}'s theorem 90.},
 fjournal = {Bulletin of the London Mathematical Society},
 journal = {Bull. Lond. Math. Soc.},
 issn = {0024-6093},
 volume = {47},
 number = {4},
 pages = {704--714},
 year = {2015},
 language = {English},
 doi = {10.1112/blms/bdv043},
 zbMATH = {6472449},
 Zbl = {1342.20053}
}

@book{Gras2003,
 author = {Gras, Georges},
 title = {Class field theory. {From} theory to practice},
 fseries = {Springer Monographs in Mathematics},
 series = {Springer Monogr. Math.},
 issn = {1439-7382},
 isbn = {3-540-44133-6},
 year = {2003},
 publisher = {Berlin: Springer},
 language = {English},
 zbMATH = {1834428},
 Zbl = {1019.11032}
}

@article{Skorobogatov2007,
 author = {Skorobogatov, A. N.},
 title = {On the elementary obstruction to the existence of rational points},
 fjournal = {Mathematical Notes},
 journal = {Math. Notes},
 issn = {0001-4346},
 volume = {81},
 number = {1},
 pages = {97--107},
 year = {2007},
 language = {English},
 doi = {10.1134/S0001434607010099},
 zbMATH = {5210297},
 Zbl = {1134.14012}
}

@article{Beauville2009,
 author = {Beauville, Arnaud},
 title = {On the {Brauer} group of {Enriques} surfaces},
 fjournal = {Mathematical Research Letters},
 journal = {Math. Res. Lett.},
 issn = {1073-2780},
 volume = {16},
 number = {5-6},
 pages = {927--934},
 year = {2009},
 language = {English},
 doi = {10.4310/MRL.2009.v16.n6.a1},
 zbMATH = {5707692},
 Zbl = {1195.14053}
}

@book{Milne1980,
 author = {Milne, J. S.},
 title = {{\'E}tale cohomology},
 fseries = {Princeton Mathematical Series},
 series = {Princeton Math. Ser.},
 volume = {33},
 year = {1980},
 publisher = {Princeton University Press, Princeton, NJ},
 language = {English},
 zbMATH = {3674235},
 Zbl = {0433.14012}
}

@book{ChaseHarrisonRosenberg1965,
 author = {Chase, S. U. and Harrison, D. K. and Rosenberg, A.},
 title = {Galois theory and cohomology of commutative rings},
 fseries = {Memoirs of the American Mathematical Society},
 series = {Mem. Am. Math. Soc.},
 issn = {0065-9266},
 volume = {52},
 isbn = {978-0-8218-1252-5; 978-0-8218-9997-7},
 year = {1965},
 publisher = {Providence, RI: American Mathematical Society (AMS)},
 language = {English},
 doi = {10.1090/memo/0052},
 zbMATH = {7073395},
 Zbl = {1415.13001}
}

@book{Neukirch2008,
 author = {Neukirch, J{\"u}rgen and Schmidt, Alexander and Wingberg, Kay},
 title = {Cohomology of number fields},
 edition = {2nd ed.},
 fseries = {Grundlehren der Mathematischen Wissenschaften},
 series = {Grundlehren Math. Wiss.},
 issn = {0072-7830},
 volume = {323},
 isbn = {978-3-540-37888-4},
 year = {2008},
 publisher = {Berlin: Springer},
 language = {English},
 zbMATH = {5162349},
 Zbl = {1136.11001}
}

@incollection{Weigel2007,
 author = {Weigel, Th.},
 title = {Frattini extensions and class field theory.},
 booktitle = {Groups St. Andrews 2005. Vol. II. Selected papers of the conference, St. Andrews, UK, July 30--August 6, 2005.},
 isbn = {978-0-521-69470-4},
 pages = {661--684},
 year = {2007},
 publisher = {Cambridge: Cambridge University Press},
 language = {English},
 zbMATH = {5155405},
 Zbl = {1119.20034}
}

@article{TorrecillasWeigel2013,
 author = {Torrecillas, B. and Weigel, Th.},
 title = {Lattices and cohomological {Mackey} functors for finite cyclic {{\(p\)}}-groups},
 fjournal = {Advances in Mathematics},
 journal = {Adv. Math.},
 issn = {0001-8708},
 volume = {244},
 pages = {533--569},
 year = {2013},
 language = {English},
 doi = {10.1016/j.aim.2013.05.010},
 zbMATH = {6264344},
 Zbl = {1284.18003}
}

@article{Aviles2020,
 author = {Gonz{\'a}lez-Avil{\'e}s, Cristian D.},
 title = {Cokernels of restriction maps and subgroups of norm one, with applications to quadratic {Galois} coverings},
 fjournal = {Journal of Pure and Applied Algebra},
 journal = {J. Pure Appl. Algebra},
 issn = {0022-4049},
 volume = {224},
 number = {2},
 pages = {559--585},
 year = {2020},
 language = {English},
 doi = {10.1016/j.jpaa.2019.05.022},
 zbMATH = {7094391},
 Zbl = {1451.14058}
}

@article{Hilbert1897,
 author = {Hilbert, D.},
 title = {Theory of algebraic number fields.},
 fjournal = {Jahresbericht der Deutschen Mathematiker-Vereinigung (DMV)},
 journal = {Jahresber. Dtsch. Math.-Ver.},
 issn = {0012-0456},
 volume = {4},
 pages = {i--xviii + 175--546},
 year = {1897},
 language = {German},
 url = {https://eudml.org/doc/144518},
 zbMATH = {2673027},
 JFM = {28.0157.05}
}

@article{Hilbert1892,
 author = {Hilbert, D.},
 title = {{\"U}ber die {Irreducibilit{\"a}t} ganzer rationaler {Functionen} mit ganzzahligen {Coefficienten}.},
 fjournal = {Journal f{\"u}r die Reine und Angewandte Mathematik},
 journal = {J. Reine Angew. Math.},
 issn = {0075-4102},
 volume = {110},
 pages = {104--129},
 year = {1892},
 language = {German},
 url = {https://eudml.org/doc/148853},
 zbMATH = {2684679},
 JFM = {24.0087.03}
}

@article{Smith2000,
 author = {Smith, Gene Ward},
 title = {Some polynomials over {{\(\mathbb{Q}(t)\)}} and their {Galois} groups},
 fjournal = {Mathematics of Computation},
 journal = {Math. Comput.},
 issn = {0025-5718},
 volume = {69},
 number = {230},
 pages = {775--796},
 year = {2000},
 language = {English},
 doi = {10.1090/S0025-5718-99-01160-6},
 zbMATH = {1414927},
 Zbl = {0962.11041}
}

@misc{Milne2020,
author={J.S. Milne},
title={Class Field Theory (v4.03)},
year={2020},
note={Available at www.jmilne.org/math/},
pages={287+viii}
}

@misc{Lang1983,
 author = {Lang, Serge},
 title = {Fundamentals of {Diophantine} geometry},
 year = {1983},
 language = {English},
 zbMATH = {3838204},
 Zbl = {0528.14013}
}

@article{Kleine2019,
 author = {Kleine, S{\"o}ren},
 title = {{{\(T\)}}-ranks of {Iwasawa} modules},
 fjournal = {Journal of Number Theory},
 journal = {J. Number Theory},
 issn = {0022-314X},
 volume = {196},
 pages = {61--86},
 year = {2019},
 language = {English},
 doi = {10.1016/j.jnt.2018.09.022},
 zbMATH = {6987925},
 Zbl = {1467.11103}
}

@article{Greenberg1973,
 author = {Greenberg, Ralph},
 title = {On a certain {{\(l\)}}-adic representation},
 fjournal = {Inventiones Mathematicae},
 journal = {Invent. Math.},
 issn = {0020-9910},
 volume = {21},
 pages = {117--124},
 year = {1973},
 language = {English},
 doi = {10.1007/BF01389691},
 url = {https://eudml.org/doc/142221},
 zbMATH = {3422497},
 Zbl = {0268.12004}
}

@article{Maksoud2023,
 author = {Maksoud, Alexandre},
 title = {On the rank of {Leopoldt}'s and {Gross}'s regulator maps},
 fjournal = {Documenta Mathematica},
 journal = {Doc. Math.},
 issn = {1431-0635},
 volume = {28},
 number = {6},
 pages = {1441--1471},
 year = {2023},
 language = {English},
 doi = {10.4171/DM/935},
 zbMATH = {7787221},
 Zbl = {1550.11112}
}

@book{Washington1997,
 author = {Washington, Lawrence C.},
 title = {Introduction to cyclotomic fields.},
 edition = {2nd ed.},
 fseries = {Graduate Texts in Mathematics},
 series = {Grad. Texts Math.},
 issn = {0072-5285},
 volume = {83},
 isbn = {0-387-94762-0},
 year = {1997},
 publisher = {New York, NY: Springer},
 language = {English},
 zbMATH = {967875},
 Zbl = {0966.11047}
}

@book{Mumford1966,
 author = {Mumford, D.},
 title = {Lectures on curves on an algebraic surface},
 fseries = {Annals of Mathematics Studies},
 series = {Ann. Math. Stud.},
 volume = {59},
 year = {1966},
 publisher = {Princeton University Press, Princeton, NJ},
 language = {English},
 doi = {10.1515/9781400882069},
 zbMATH = {3299277},
 Zbl = {0187.42701}
}

@misc{OSCAR,
  key          = {OSCAR},
  organization = {The OSCAR Team},
  title        = {{OSCAR} -- {O}pen {S}ource {C}omputer {A}lgebra {R}esearch system, {V}ersion 1.5.0-DEV},
  year         = {2025},
  url          = {https://www.oscar-system.org},
}

@article{GrasMunnier1998,
     author = {Georges Gras and Adeline Munnier},
     title = {Extensions cycliques $T$-totalement ramifi\'ees},
     journal = {Publications math\'ematiques de Besan\c{c}on. Alg\`ebre et th\'eorie des nombres},
     eid = {6},
     pages = {1--17},
     publisher = {Presses universitaires de Franche-Comt\'e},
     year = {1998},
     doi = {10.5802/pmb.a-91},
     language = {fr},
     url = {https://pmb.centre-mersenne.org/articles/10.5802/pmb.a-91/}
}

@article{Kuzmin1972,
 author = {Kuz'min, L. V.},
 title = {The {Tate} module for algebraic number fields},
 fjournal = {Izvestiya Akademii Nauk SSSR. Seriya Matematicheskaya},
 journal = {Izv. Akad. Nauk SSSR, Ser. Mat.},
 issn = {0373-2436},
 volume = {36},
 pages = {267--327},
 year = {1972},
 language = {Russian},
 zbMATH = {3365393},
 Zbl = {0231.12013}
}

@article{Jaulent1994,
 author = {Jaulent, Jean-Fran{\c{c}}ois},
 title = {Logarithmic classes of number fields},
 fjournal = {Journal de Th{\'e}orie des Nombres de Bordeaux},
 journal = {J. Th{\'e}or. Nombres Bordx.},
 issn = {1246-7405},
 volume = {6},
 number = {2},
 pages = {301--325},
 year = {1994},
 language = {French},
 doi = {10.5802/jtnb.117},
 url = {https://eudml.org/doc/93606},
 zbMATH = {754095},
 Zbl = {0827.11064}
}

@misc{Jaulent1986,
 author = {Jaulent, Jean-Fran{\c{c}}ois},
 title = {L'arithm{\'e}tique des {{\(\ell\)}}-extensions.},
 year = {1986},
 language = {French},
 zbMATH = {3968700},
 Zbl = {0601.12002}
}

@article{Gross1981,
 author = {Gross, Benedict H.},
 title = {p-adic {L}-series at s=0},
 fjournal = {Journal of the Faculty of Science. Section I A},
 journal = {J. Fac. Sci., Univ. Tokyo, Sect. I A},
 issn = {0040-8980},
 volume = {28},
 pages = {979--994},
 year = {1981},
 language = {English},
 zbMATH = {3799778},
 Zbl = {0507.12010}
}

@article{Fukuda1994,
 author = {Fukuda, Takashi},
 title = {Remarks on {{\(\mathbb{Z}_ p\)}}-extensions of number fields},
 fjournal = {Proceedings of the Japan Academy. Series A},
 journal = {Proc. Japan Acad., Ser. A},
 issn = {0386-2194},
 volume = {70},
 number = {8},
 pages = {264--266},
 year = {1994},
 language = {English},
 doi = {10.3792/pjaa.70.264},
 zbMATH = {723970},
 Zbl = {0823.11064}
}

@book{Cohen2000,
 author = {Cohen, Henri},
 title = {Advanced topics in computational number theory},
 fseries = {Graduate Texts in Mathematics},
 series = {Grad. Texts Math.},
 issn = {0072-5285},
 volume = {193},
 isbn = {0-387-98727-4},
 year = {2000},
 publisher = {New York, NY: Springer},
 language = {English},
 doi = {10.1007/978-1-4419-8489-0},
 zbMATH = {1391020},
 Zbl = {0977.11056}
}

@article{Jaulent2017,
 author = {Jaulent, Jean-Fran{\c{c}}ois},
 title = {On the cyclotomic norms and the {Leopoldt} and {Gross}-{Kuz}'min conjectures},
 fjournal = {Annales Math{\'e}matiques du Qu{\'e}bec},
 journal = {Ann. Math. Qu{\'e}.},
 issn = {2195-4755},
 volume = {41},
 number = {1},
 pages = {119--140},
 year = {2017},
 language = {French},
 doi = {10.1007/s40316-016-0069-3},
 zbMATH = {6760539},
 Zbl = {1432.11161}
}

@article{Maksoud2024,
 author = {Maksoud, Alexandre},
 title = {On generalized main conjectures and {{\(p\)}}-adic {Stark} conjectures for {Artin} motives},
 fjournal = {Transactions of the American Mathematical Society},
 journal = {Trans. Am. Math. Soc.},
 issn = {0002-9947},
 volume = {377},
 number = {6},
 pages = {3845--3904},
 year = {2024},
 language = {English},
 doi = {10.1090/tran/9131},
 zbMATH = {7853745},
 Zbl = {1551.11120}
}

@article{Jaulent2002,
 author = {Jaulent, Jean-Fran{\c{c}}ois},
 title = {Logarithmic classes of totally real number fields},
 fjournal = {Acta Arithmetica},
 journal = {Acta Arith.},
 issn = {0065-1036},
 volume = {103},
 number = {1},
 pages = {1--7},
 year = {2002},
 language = {French},
 doi = {10.4064/aa103-1-1},
 zbMATH = {1789546},
 Zbl = {1001.11044}
}

@article{Kuzmin2018,
 author = {Kuz'min, L. V.},
 title = {{Local and global universal norms in the cyclotomic \( \mathbb{Z}_{\ell}\)-extension of an algebraic number field}},
 fjournal = {Izvestiya: Mathematics},
 journal = {Izv. Math.},
 issn = {1064-5632},
 volume = {82},
 number = {3},
 pages = {532--548},
 year = {2018},
 language = {English},
 doi = {10.1070/IM8653},
 zbMATH = {6902317},
 Zbl = {1439.11277}
}

@misc{Gras2024,
 author = {Gras, Georges},
 title = {{On the $\mathbb{Z}_p$-extensions of a totally $p$-adic imaginary quadratic field -- {With} an appendix by {Jean}-{Fran{\c{c}}ois} {Jaulent}}},
 year = {2024},
 howpublished = {Preprint, {arXiv}:2403.16603 [math.{NT}] (2024)},
 url = {https://arxiv.org/abs/2403.16603},
 arXiv = {arXiv:2403.16603}
}

@article{Aviles2004,
 author = {Gonzalez-Avil{\'e}s, Cristian D.},
 title = {Finite modules over non-semisimple group rings.},
 fjournal = {Israel Journal of Mathematics},
 journal = {Isr. J. Math.},
 issn = {0021-2172},
 volume = {144},
 pages = {61--92},
 year = {2004},
 language = {English},
 doi = {10.1007/BF02984406},
 zbMATH = {2147122},
 Zbl = {1086.20005}
}

@misc{Dress1973,
 author = {Dress, Andreas W. M.},
 title = {Contributions to the theory of induced representations},
 year = {1973},
 language = {English},
 howpublished = {Algebr. {{\(K\)}}-{Theory} {II}, {Proc}. {Conf}. {Battelle} {Inst}. 1972, {Lect}. {Notes} {Math}. 342, 183-240 (1973).},
 zbMATH = {3517377},
 Zbl = {0331.18016}
}

@article{Gras2016,
 author = {Gras, Georges},
 title = {The local {{\(\theta\)}}-regulators of an algebraic number: {{\(p\)}}-adic conjectures},
 fjournal = {Canadian Journal of Mathematics},
 journal = {Can. J. Math.},
 issn = {0008-414X},
 volume = {68},
 number = {3},
 pages = {571--624},
 year = {2016},
 language = {French},
 doi = {10.4153/CJM-2015-026-3},
 zbMATH = {6589334},
 Zbl = {1351.11033}
}

@article{Movahhedi1990,
 author = {Movahhedi, A.},
 title = {Sur les p-extensions des corps p-rationnels. ({On} p-extensions of p- rational fields)},
 fjournal = {Mathematische Nachrichten},
 journal = {Math. Nachr.},
 issn = {0025-584X},
 volume = {149},
 pages = {163--176},
 year = {1990},
 language = {French},
 doi = {10.1002/mana.19901490113},
 zbMATH = {4191814},
 Zbl = {0723.11054}
}

@phdthesis{Bembom2012,
 author = {Bembom, T.},
 title = {The capitulation problem in class field theory},
 year = {2012},
 publisher = {G{\"o}ttingen: Univ. G{\"o}ttingen (Diss.)},
 language = {English},
 school ={Georg-August-Universität Göttingen},
 url = {ediss.uni-goettingen.de/bitstream/handle/11858/00-1735-0000-000D-F05F-8/bembom.pdf?sequence=1},
 zbMATH = {6346838},
 Zbl = {1298.11104}
}

@misc{Feuerpfeil2026Code,
author = {Feuerpfeil, J.},
title = {{ Hilbert 90 for S Class Groups}},
year = {2026},
note = {https://github.com/JulianFeuerpfeil/Hilbert-90-for-S-Class-Groups}
}

@article{BleyBoltje2004,
 author = {Bley, Werner and Boltje, Robert},
 title = {Cohomological {Mackey} functors in number theory},
 fjournal = {Journal of Number Theory},
 journal = {J. Number Theory},
 issn = {0022-314X},
 volume = {105},
 number = {1},
 pages = {1--37},
 year = {2004},
 language = {English},
 doi = {10.1016/j.jnt.2003.09.002},
 zbMATH = {2082492},
 Zbl = {1061.11059}
}

\end{document}